\input amstex
\documentstyle{amsppt}
%
\catcode`@=11
\redefine\output@{%
  \def\break{\penalty-\@M}\let\par\endgraf
  \ifodd\pageno\global\hoffset=105pt\else\global\hoffset=8pt\fi  
  \shipout\vbox{%
    \ifplain@
      \let\makeheadline\relax \let\makefootline\relax
    \else
      \iffirstpage@ \global\firstpage@false
        \let\rightheadline\frheadline
        \let\leftheadline\flheadline
      \else
        \ifrunheads@ 
        \else \let\makeheadline\relax
        \fi
      \fi
    \fi
    \makeheadline \pagebody \makefootline}%
  \advancepageno \ifnum\outputpenalty>-\@MM\else\dosupereject\fi
}
\catcode`\@=\active
\nopagenumbers
\def\negskp{\hskip -2pt}
\def\MatGrU{\operatorname{U}}
\def\MatGrSU{\operatorname{SU}}
\def\MatAlgSU{\operatorname{su}}
\def\Alpha{\operatorname{A}}
\def\idop{\operatorname{\bold{id}}}
\def\const{\operatorname{const}}
\def\msum#1{\operatornamewithlimits{\sum^#1\!{\ssize\ldots}\!\sum^#1}}
\def\compos{\,\raise 1pt\hbox{$\sssize\circ$} \,}
\accentedsymbol\bulletH{\overset{\kern 2pt\sssize\bullet}\to H}
\accentedsymbol\circH{\overset{\kern 2pt\sssize\circ}\to H}
\accentedsymbol\bulletBH{\overset{\sssize\bullet}\to{\bold H}}
\accentedsymbol\circBH{\overset{\sssize\circ}\to{\bold H}}
\loadbold
\accentedsymbol\bulpsi{\overset{\kern 2pt\sssize\bullet\kern -2pt}%
   \to\psi}
\accentedsymbol\circpsi{\overset{\kern 2pt\sssize\circ\kern -2pt}%
   \to\psi}
\def\vtrule{\vrule height 12pt depth 6pt}
\def\chirk{\special{em:point 1}\kern 1.2pt\raise 0.6pt
  \hbox to 0pt{\special{em:point 2}\hss}\kern -1.2pt
  \special{em:line 1,2,0.3pt}\ignorespaces}
\def\Chirk{\special{em:point 1}\kern 1.5pt\raise 0.6pt
  \hbox to 0pt{\special{em:point 2}\hss}\kern -1.5pt
  \special{em:line 1,2,0.3pt}\ignorespaces}
\def\Hirk{\kern 0pt\special{em:point 1}\kern 4pt\special{em:point 2}
  \kern -3pt\special{em:line 1,2,0.3pt}\ignorespaces}
\accentedsymbol\uuud{d\unskip\kern -3.8pt\raise 1pt\hbox to
  0pt{\chirk\hss}\kern 0.2pt\raise 1.9pt\hbox to 0pt{\chirk\hss}
  \kern 0.2pt\raise 2.8pt\hbox to 0pt{\chirk\hss}\kern 3.5pt}
\accentedsymbol\bolduuud{\bold d\unskip\kern -4pt\raise 1.1pt\hbox to
  0pt{\chirk\hss}\kern 0pt\raise 1.9pt\hbox to 0pt{\chirk\hss}
  \kern 0pt\raise 2.5pt\hbox to 0pt{\chirk\hss}\kern 4pt}
\accentedsymbol\uud{d\unskip\kern -3.72pt\raise 1.5pt\hbox to
  0pt{\chirk\hss}\kern 0.26pt\raise 2.4pt\hbox to 0pt{\chirk\hss}
  \kern 3.5pt}
\accentedsymbol\bolduud{\bold d\unskip\kern -4pt\raise 1.4pt\hbox to
  0pt{\chirk\hss}\kern 0pt\raise 2.3pt\hbox to 0pt{\chirk\hss}
  \kern 4pt}
\accentedsymbol\uuuD{D\unskip\kern -4.6pt\raise 2pt\hbox to
  0pt{\Chirk\hss}\kern 0.2pt\raise 3pt\hbox to 0pt{\Chirk\hss}
  \kern 0.2pt\raise 4.1pt\hbox to 0pt{\Chirk\hss}\kern 4.5pt}
\accentedsymbol\uuD{D\unskip\kern -4.5pt\raise 2.4pt\hbox to
  0pt{\Chirk\hss}\kern 0.2pt\raise 3.4pt\hbox to 0pt{\Chirk\hss}
  \kern 4.4pt}
\accentedsymbol\uD{D\unskip\kern -4.4pt\raise 3pt\hbox to
  0pt{\Chirk\hss}\kern 4.4pt}
\accentedsymbol\bolduuuD{\bold D\unskip\kern -4.6pt\raise 2pt\hbox to
  0pt{\Chirk\hss}\kern 0pt\raise 3pt\hbox to 0pt{\Chirk\hss}
  \kern 0pt\raise 4.1pt\hbox to 0pt{\Chirk\hss}\kern 4.5pt}
\accentedsymbol\bolduuD{\bold D\unskip\kern -4.5pt\raise 2.4pt\hbox to
  0pt{\Chirk\hss}\kern 0pt\raise 3.4pt\hbox to 0pt{\Chirk\hss}
  \kern 4.4pt}
\accentedsymbol\bolduD{\bold D\unskip\kern -4.5pt\raise 3pt\hbox to
  0pt{\Chirk\hss}\kern 4.4pt}
\accentedsymbol\uuuU{U\unskip\kern -5.1pt\raise 2pt\hbox to
  0pt{\Chirk\hss}\kern 0.2pt\raise 3pt\hbox to 0pt{\Chirk\hss}
  \kern 0.2pt\raise 4.1pt\hbox to 0pt{\Chirk\hss}\kern 4.9pt}
\accentedsymbol\uuU{U\unskip\kern -5.0pt\raise 2.4pt\hbox to
  0pt{\Chirk\hss}\kern 0.2pt\raise 3.4pt\hbox to 0pt{\Chirk\hss}
  \kern 4.9pt}
\accentedsymbol\uU{U\unskip\kern -4.9pt\raise 3pt\hbox to
  0pt{\Chirk\hss}\kern 4.8pt}
\accentedsymbol\bolduPsi{\boldsymbol\Psi\unskip\kern -6.5pt\raise 3pt
  \hbox to 0pt{\Hirk\hss}\kern 6.6pt}
\accentedsymbol\bolduuPsi{\boldsymbol\Psi\unskip\kern -6.5pt\raise 3pt
  \hbox to 0pt{\Hirk\hss}\unskip\kern 0pt\raise 3.8pt
  \hbox to 0pt{\Hirk\hss}\kern 6.6pt}
\accentedsymbol\bolduuuPsi{\boldsymbol\Psi\unskip\kern -6.5pt\raise 3pt
  \hbox to 0pt{\Hirk\hss}\unskip\kern 0pt\raise 3.8pt
  \hbox to 0pt{\Hirk\hss}\unskip\kern 0pt\raise 4.6pt
  \hbox to 0pt{\Hirk\hss}\kern 6.6pt}
\accentedsymbol\uA{\operatorname{A}\unskip\kern -6.2pt\raise 2.8pt\hbox to
  0pt{\vbox{\hrule width 2pt height 0.3pt}\hss}\kern 4.5pt}
\accentedsymbol\uuA{\operatorname{A}\unskip\kern -6.2pt\raise 2.8pt\hbox 
  to 0pt{\vbox{\hrule width 2pt height 0.3pt}\hss}\unskip\kern 0.2pt\raise
  3.4pt\hbox to 0pt{\vbox{\hrule width 1.6pt height 0.3pt}\hss}\kern 4.3pt}
\accentedsymbol\uuuA{\operatorname{A}\unskip\kern -6.2pt\raise 2.8pt\hbox 
  to 0pt{\vbox{\hrule width 2pt height 0.3pt}\hss}\unskip\kern 0.2pt\raise
  3.4pt\hbox to 0pt{\vbox{\hrule width 1.6pt height 0.3pt}\hss}\unskip
  \kern 0.2pt\raise 4.0pt\hbox to 0pt{\vbox{\hrule width 1.0pt height
  0.3pt}\hss}\kern 4.3pt}
\accentedsymbol\uF{\Cal F\unskip\kern -6.4pt\raise
  3.4pt\hbox to 0pt{\vbox{\hrule width 1.6pt height 0.3pt}\hss}\kern 6.7pt}
\accentedsymbol\uuF{\Cal F\unskip\kern -6.6pt\raise 2.8pt\hbox 
  to 0pt{\vbox{\hrule width 2pt height 0.3pt}\hss}\unskip\kern 0.2pt\raise
  3.4pt\hbox to 0pt{\vbox{\hrule width 1.6pt height 0.3pt}\hss}\kern 6.9pt}
\accentedsymbol\uuuF{\Cal F\unskip\kern -6.6pt\raise 2.8pt\hbox 
  to 0pt{\vbox{\hrule width 2pt height 0.3pt}\hss}\unskip\kern 0.2pt\raise
  3.4pt\hbox to 0pt{\vbox{\hrule width 1.6pt height 0.3pt}\hss}\unskip
  \kern 0.2pt\raise 4.0pt\hbox to 0pt{\vbox{\hrule width 1.45pt height
  0.3pt}\hss}\kern 6.9pt}
\accentedsymbol\uCA{\Cal A\unskip\kern -5.6pt\raise
  3.4pt\hbox to 0pt{\vbox{\hrule width 1.6pt height 0.3pt}\hss}\kern 5.9pt}
\accentedsymbol\uuCA{\Cal A\unskip\kern -5.9pt\raise 2.8pt\hbox 
  to 0pt{\vbox{\hrule width 2pt height 0.3pt}\hss}\unskip\kern 0.2pt\raise
  3.4pt\hbox to 0pt{\vbox{\hrule width 1.6pt height 0.3pt}\hss}\kern 6.2pt}
\accentedsymbol\uuuCA{\Cal A\unskip\kern -5.9pt\raise 2.8pt\hbox 
  to 0pt{\vbox{\hrule width 2pt height 0.3pt}\hss}\unskip\kern 0.3pt\raise
  3.4pt\hbox to 0pt{\vbox{\hrule width 1.6pt height 0.3pt}\hss}\unskip
  \kern 0.5pt\raise 4.0pt\hbox to 0pt{\vbox{\hrule width 1.45pt height
  0.3pt}\hss}\kern 6.2pt}
\loadeurb
\accentedsymbol\bolduF{\eurb F\unskip\kern -5.5pt\raise
  2.8pt\hbox to 0pt{\vbox{\hrule width 1.6pt height 0.3pt}\hss}\kern 5.3pt}
\accentedsymbol\bolduuF{\eurb F\unskip\kern -5.5pt\raise 2.8pt\hbox 
  to 0pt{\vbox{\hrule width 1.6pt height 0.3pt}\hss}\unskip\kern 0pt\raise 
  3.4pt\hbox to 0pt{\vbox{\hrule width 1.6pt height 0.3pt}\hss}\kern 5.3pt}
\accentedsymbol\bolduuuF{\eurb F\unskip\kern -5.5pt\raise 2.8pt\hbox 
  to 0pt{\vbox{\hrule width 1.6pt height 0.3pt}\hss}\unskip\kern 0pt\raise 
  3.4pt\hbox to 0pt{\vbox{\hrule width 1.6pt height 0.3pt}\hss}\unskip
  \kern 0pt\raise 4.0pt\hbox to 0pt{\vbox{\hrule width 1.6pt height
  0.3pt}\hss}\kern 5.3pt}
\accentedsymbol\bolduCA{\eurb A\unskip\kern -7.9pt\raise
  2.8pt\hbox to 0pt{\vbox{\hrule width 1.6pt height 0.3pt}\hss}\kern 7.7pt}
\accentedsymbol\bolduuCA{\eurb A\unskip\kern -7.9pt\raise 2.8pt\hbox 
  to 0pt{\vbox{\hrule width 1.6pt height 0.3pt}\hss}\unskip\kern 0.3pt
  \raise 3.4pt\hbox to 0pt{\vbox{\hrule width 1.6pt height 0.3pt}\hss}
  \kern 7.7pt}
\accentedsymbol\bolduuuCA{\eurb A\unskip\kern -7.9pt\raise 2.8pt\hbox 
  to 0pt{\vbox{\hrule width 1.6pt height 0.3pt}\hss}\unskip\kern 0.3pt
  \raise 3.4pt\hbox to 0pt{\vbox{\hrule width 1.6pt height 0.3pt}\hss}
  \unskip\kern 0.3pt\raise 4.0pt\hbox to 0pt{\vbox{\hrule width 1.6pt 
  height 0.3pt}\hss}\kern 7.7pt}
\def\blue#1{#1}
\catcode`#=11\def\diez{#}\catcode`#=6
\catcode`_=11\def\podcherkivanie{_}\catcode`_=8
\def\mycite#1{\cite{\blue{#1}}\immediate\special{ps:
     ShrHPSdict begin /ShrBORDERthickness 0 def}}

\def\mytag#1{%
    \tag#1}
\def\mythetag#1{\thetag{\blue{#1}}\immediate\special{ps:
     ShrHPSdict begin /ShrBORDERthickness 0 def}}
\def\myrefno#1{\no#1}
\def\myhref#1#2{\blue{#2}\immediate\special{ps:
     ShrHPSdict begin /ShrBORDERthickness 0 def}}
\def\myEarXivlink{\myhref{http://arXiv.org}{http:/\negskp/arXiv.org}}
\def\myGeoCities{\myhref{http://www.geocities.com}{GeoCities}}
\def\mytheorem#1{\csname proclaim\endcsname{Theorem #1}}

\def\mylemma#1{\csname proclaim\endcsname{Lemma #1}}

\def\mycorollary#1{\csname proclaim\endcsname{Corollary #1}}

\def\mydefinition#1{\definition{Definition #1}}
\def\mythedefinition#1{\blue{#1}\immediate\special{ps:
     ShrHPSdict begin /ShrBORDERthickness 0 def}}

\pagewidth{360pt}
\pageheight{606pt}
\topmatter
\title
A note on the Standard Model\\in a gravitation field.
\endtitle
\author
R.~A.~Sharipov
\endauthor
\address 5 Rabochaya street, 450003 Ufa, Russia\newline
\vphantom{a}\kern 12pt Cell Phone: +7-(917)-476-93-48
\endaddress
\email \vtop to 30pt{\hsize=280pt\noindent
\myhref{mailto:r-sharipov\@mail.ru}
{r-sharipov\@mail.ru}\newline
\myhref{mailto:ra\podcherkivanie sharipov\@lycos.com}
{ra\_\hskip 1pt sharipov\@lycos.com}\newline
\myhref{mailto:R\podcherkivanie Sharipov\@ic.bashedu.ru}
{R\_\hskip 1pt Sharipov\@ic.bashedu.ru}\vss}
\endemail
\urladdr
\vtop to 20pt{\hsize=280pt\noindent
\myhref{http://www.geocities.com/r-sharipov}
{http:/\negskp/www.geocities.com/r-sharipov}\newline
\myhref{http://www.freetextbooks.boom.ru/index.html}
{http:/\negskp/www.freetextbooks.boom.ru/index.html}\vss}
\endurladdr
\abstract
    The Standard Model of elementary particles is a theory unifying
three of the four basic forces of the Nature: electromagnetic, weak,
and strong interactions. In this paper we consider the Standard
Model in the presence of a classical (non-quantized) gravitation field
and apply a bundle approach for describing it.
\endabstract
\endtopmatter
\loadeufb
\TagsOnRight
\document

\rightheadtext{A note on the Standard Model \dots}
\head
1. Fermion fields of the Standard Model. 
\endhead
    Fermion fields of the Standard Model are subdivided into two parts:
lepton fields and quark fields. Lepton fields are subdivided into three
generations. The first generation is represented by an electron $e$ and 
an electronic neutrino $\nu_e$, the second generation is represented by 
a muon $\mu$ and its neutrino $\nu_\mu$, and the third generation is
represented by a tauon $\tau$ and its neutrino $\nu_\tau$.
$$
\vcenter{\hsize 10cm
\offinterlineskip\settabs\+\indent
\vtrule
\hskip 3.2cm &\vtrule 
\hskip 3.2cm &\vtrule 
\hskip 3.2cm &\vtrule
\cr\hrule 
\+\vtrule
\hfill 1-st generation\hfill&\vtrule
\hfill 2-nd generation\hfill &\vtrule
\hfill 3-rd generation\hfill &\vtrule
\cr\hrule 
\+\vtrule
\hfill $e$-neutrino $\nu_e$\hfill&\vtrule
\hfill $\mu$-neutrino $\nu_\mu$\hfill &\vtrule
\hfill $\tau$-neutrino $\nu_\tau$\hfill &\vtrule
\cr\hrule 
\+\vtrule
\hfill electron $e$\hfill&\vtrule
\hfill muon $\mu$\hfill &\vtrule
\hfill tauon $\tau$\hfill &\vtrule
\cr\hrule 
}\quad
\mytag{1.1}
$$
In a similar way, quarks are subdivided into three generations. They 
are represented in the following table similar to the above table
\mythetag{1.1}:
$$
\vcenter{\hsize 10cm
\offinterlineskip\settabs\+\indent
\vtrule
\hskip 3.2cm &\vtrule 
\hskip 3.2cm &\vtrule 
\hskip 3.2cm &\vtrule
\cr\hrule 
\+\vtrule
\hfill 1-st generation\hfill&\vtrule
\hfill 2-nd generation\hfill &\vtrule
\hfill 3-rd generation\hfill &\vtrule
\cr\hrule 
\+\vtrule
\hfill up-quark $u$\hfill&\vtrule
\hfill charm-quark $c$\hfill &\vtrule
\hfill top-quark $t$\hfill &\vtrule
\cr\hrule 
\+\vtrule
\hfill down-quark $d$\hfill&\vtrule
\hfill strange-quark $s$\hfill &\vtrule
\hfill bottom-quark $b$\hfill &\vtrule
\cr\hrule 
}\quad
\mytag{1.2}
$$
Leptons participate in electromagnetic and weak interactions. These
interactions are described by the $\MatGrU(1)\times\MatGrSU(2)$ symmetry
which is spontaneously broken according to the Higgs mechanism.
Moreover, they break the chiral symmetry on the level of Dirac spinors.
This symmetry is often called the left-to-right symmetry, but we prefer
to say the chiral-to-antichiral symmetry or simply the chiral symmetry
(see some details in \mycite{1}). We distinguish between lepton wave
functions by mens of the generation index enclosed into square brackets:
$$
\xalignat 3
&\hskip -2em
\boldsymbol\psi[e],
&&\boldsymbol\psi[\mu],
&&\boldsymbol\psi[\tau].
\mytag{1.3}
\endxalignat
$$
The wave functions \mythetag{1.3} have chiral and antichiral constituent
parts. Chiral parts are doublets with respect to $\MatGrSU(2)$ symmetry:
$$
\xalignat 3
&\hskip -2em
\bulpsi^{a\alpha}_{111}[e],
&&\bulpsi^{a\alpha}_{111}[\mu],
&&\bulpsi^{a\alpha}_{111}[\tau].
\mytag{1.4}
\endxalignat
$$
The indices in \mythetag{1.4} means that we take the space-time manifold
$M$ equipped the appropriate complex vector bundles $DM$, $\uU\!M$, 
and $S\uuU\!M$ (see \mycite{2}). The index $a=1,2,3,4$ in \mythetag{1.4} 
is a spinor index associated with the Dirac bundle $DM$ or, more
precisely, with some frame $(U,\,\boldsymbol\Psi_1,\,\boldsymbol\Psi_2,
\,\boldsymbol\Psi_3,\,\boldsymbol\Psi_4)$ of $DM$. The index $\alpha=1,2$
is associated with some frame $(U,\,\bolduuPsi_1,\,\bolduuPsi_2)$ of the
two-dimensional complex bundle $S\uuU\!M$. Due to the presence of this
index the wave functions \mythetag{1.4} are said to be
$\MatGrSU(2)$-doublets. Three lower indices in \mythetag{1.4} are always
equal to unity because they are associated with some frame $(U,\,
\bolduPsi_1)$ of the one-dimensional bundle $\uU\!M$.\par
     The antichiral parts of the wave functions \mythetag{1.3} are
$\MatGrSU(2)$-singlets. Their components have one spinor index $a$
and six $\MatGrU(1)$ indices equal to $1$:
$$
\xalignat 3
&\hskip -2em
\circpsi^a_{111111}[e],
&&\circpsi^a_{111111}[\mu],
&&\circpsi^a_{111111}[\tau].
\mytag{1.5}
\endxalignat
$$
By usual convention (see \mycite{3}) antichiral (right) neutrinos are not
considered. In this paper we follow this convention, though in some papers
right neutrinos are introduced, e\.\,g\. in \mycite{4}.\par
     The wave functions \mythetag{1.4} and \mythetag{1.5} are chiral and
antichiral in the sense of the following equalities relating them
with the components of the chirality operator $\bold H$:
$$
\xalignat 2
&\hskip -2em
\sum^4_{a=1}H^b_a\,\bulpsi^{a\alpha}_{111}[q]
=\bulpsi^{b\alpha}_{111}[e],
&&\sum^4_{a=1}H^b_a\,\circpsi^a_{111111}[q]
=-\circpsi^b_{111111}[q].\quad
\mytag{1.6}
\endxalignat
$$
Here $q=e,\mu,\tau$ is the generation index. The chirality operator
$\bold H$ is an attribute of the Dirac bundle (see details in \mycite{1}).
In physical literature (see \mycite{3} as an example), when the flat
Minkowski space is taken for the space-time manifold $M$, the chirality
operator $\bold H$ is represented by the Dirac matrix $\gamma^5$
$$
\hskip -2em
H^b_a=\Vmatrix 1 & 0 & 0 & 0\\0 & 1 & 0 & 0\\
0 & 0 & -1 & 0\\0 & 0 & 0 & -1\endVmatrix
=\gamma^5=i\,\gamma^0\,\gamma^1\,\gamma^2\,\gamma^3.
\mytag{1.7}
$$
Other Dirac $\gamma$-matrices are given by the following formulas
(see \thetag{1.13} in \mycite{5}):
$$
\xalignat 2
&\hskip -2em
\gamma^{b\,0}_a=\Vmatrix 0&0&1&0\\0&0&0&1\\1&0&0&0\\0&1&0&0\endVmatrix,
&&\gamma^{b1}_a=\Vmatrix 0&0&0&-1\\0&0&-1&0\\0&1&0&0\\1&0&0&0\endVmatrix,
\quad\\
\vspace{-1.5ex}
&&&\mytag{1.8}\\
\vspace{-1.5ex}
&\hskip -2em
\gamma^{b\,2}_a=\Vmatrix 0&0&0&i\\0&0&-i&0\\0&-i&0&0\\i&0&0&0\endVmatrix,
&&\gamma^{b\,3}_a=\Vmatrix 0&0&-1&0\\0&0&0&1\\1&0&0&0\\0&-1&0&0\endVmatrix.
\quad
\endxalignat
$$
By means of the chirality operator $\bold H$ we define two projection
operators:
$$
\xalignat 2
&\hskip -2em
\bulletBH=\frac{\idop+\bold H}{2},
&&\circBH=\frac{\idop-\bold H}{2}.
\mytag{1.9}
\endxalignat
$$
Here $\idop$ is the identity operator. Therefore, the components of 
the projection operators \mythetag{1.9} are given by the formulas
$$
\xalignat 2
&\hskip -2em
\bulletH^b_a=\frac{\delta^b_a+H^b_a}{2},
&&\circH^b_a=\frac{\delta^b_a-H^b_a}{2}.
\mytag{1.10}
\endxalignat
$$
By means of \mythetag{1.10} the equalities \mythetag{1.6} are written
as follows:
$$
\xalignat 2
&\hskip -2em
\sum^4_{a=1}\circH^b_a\,\bulpsi^{a\alpha}_{111}[q]=0,
&&\sum^4_{a=1}\bulletH^b_a\,\circpsi^a_{111111}[q]=0.\quad
\mytag{1.11}
\endxalignat
$$
The indices $a$ and $b$ in \mythetag{1.6}, \mythetag{1.7}, \mythetag{1.8},
\mythetag{1.9}, \mythetag{1.10}, and \mythetag{1.11} are spinor indices.
They are associated with the Dirac bundle $DM$. The third index of the 
Dirac matrices represented by the numbers $0,1,2,3$ in \mythetag{1.8} is
a spacial index, it is associated with the tangent bundle $TM$.\par
    Like in the case of leptons, quark wave functions are subdivided into
three generations according to the generation table \mythetag{1.2}:
$$
\xalignat 3
&\hskip -2em
\boldsymbol\psi[1],
&&\boldsymbol\psi[2],
&&\boldsymbol\psi[3].\quad
\mytag{1.12}
\endxalignat
$$
However, now we use a numeric index for generations, since
$\boldsymbol\psi[1]$ describes both an up-quark and a down-quark. 
Similarly, $\boldsymbol\psi[2]$ describes a charm-quark together with 
a strange-quark and $\boldsymbol\psi[3]$ describes a top-quark together
with a bottom-quark. Chiral and antichiral parts of the wave functions 
\mythetag{1.12} behave differently with respect to the $\MatGrSU(2)$
symmetry. Chiral parts are $\MatGrSU(2)$-doublets:
$$
\xalignat 3
&\hskip -2em
\bulpsi^{a1\alpha\beta}[1],
&&\bulpsi^{a1\alpha\beta}[2],
&&\bulpsi^{a1\alpha\beta}[3].
\mytag{1.13}
\endxalignat
$$
Note that in \mythetag{1.13}, in contrast to \mythetag{1.4}, we have one
more index. The additional index $\beta=1,2,3$ is responsible for color,
it describes strong interactions of quarks. For antichiral parts of the
wave functions \mythetag{1.12} the index $\alpha$ is omitted:
$$
\xalignat 3
&\hskip -2em
\circpsi^{a1111\beta}[u],
&&\circpsi^{a1111\beta}[c],
&&\circpsi^{a1111\beta}[t],\\
\vspace{-1.5ex}
&&&\mytag{1.14}\\
\vspace{-1.5ex}
&\hskip -2em
\circpsi^{a\beta}_{11}[d],
&&\circpsi^{a\beta}_{11}[s],
&&\circpsi^{a\beta}_{11}[b].
\endxalignat
$$
They are $\MatGrSU(2)$-singlets. The wave functions
\mythetag{1.13} and \mythetag{1.14} are chiral and antichiral in the 
sense of the following equalities:
$$
\xalignat 2
&\hskip -2em
\sum^4_{a=1}H^b_a\,\bulpsi^{a1\alpha\beta}[q]=\bulpsi^{b\kern 0.4pt 1
\alpha\beta}[q],
&&q=1,2,3;\\
&\hskip -2em
\sum^4_{a=1}H^b_a\,\circpsi^{a1111\beta}[q]=-\circpsi^{b\kern 0.4pt
1111\beta}[q],&&q=u,c,t;
\mytag{1.15}\\
&\hskip -2em
\sum^4_{a=1}H^b_a\,\circpsi^{a\beta}_{11}[q]=-\circpsi^{b\kern 0.4pt
\beta}_{11}[q],&&q=d,s,b.
\endxalignat
$$
The equalities \mythetag{1.15} are analogous to \mythetag{1.6}. In terms 
of the projection operators introduced by the formulas \mythetag{1.9} they
are rewritten as
$$
\xalignat 3
&\hskip -1em
\sum^4_{a=1}\circH^b_a\,\bulpsi^{a1\alpha\beta}[q]=0,
&&\sum^4_{a=1}\bulletH^b_a\,\circpsi^{a1111\beta}[q]=0,
&&\sum^4_{a=1}\bulletH^b_a\,\circpsi^{a\beta}_{11}[q]=0.
\qquad\qquad
\mytag{1.16}
\endxalignat
$$
In this form \mythetag{1.16} the identities \mythetag{1.15} are
analogous to the identities \mythetag{1.11}.\par
\head
2. The Higgs field and the classical vacuum. 
\endhead
     The classical Higgs field $\boldsymbol\varphi$ is a scalar field 
being $\MatGrSU(2)$-doublet with respect to electro-weak bundle $S\uuU\!M$.
In a coordinate form, i\.\,e\. upon choosing some frame $(U,\,\bolduuPsi_1,
\,\bolduuPsi_2)$ of $S\uuU\!M$ and some frame $(U,\,\bolduPsi_1)$ of
$\uU\!M$, it is represented as 
$$
\hskip -2em
\varphi^{\kern 0.4pt\alpha111}\text{, \ where \ }\alpha=1,2. 
\mytag{2.1}
$$
Let $(U,\,\boldsymbol\Upsilon_0,\,\boldsymbol\Upsilon_1,\,\boldsymbol
\Upsilon_2,\,\boldsymbol\Upsilon_3)$ be some frame of the tangent bundle
$TM$ and let $L_{\boldsymbol\Upsilon_k}$ be the Lie derivative along
the vector field $\boldsymbol\Upsilon_k$. Then we can define a covariant
differentiation of the Higgs field $\boldsymbol\varphi$ through the
following formula:
$$
\hskip -2em
\nabla_{\!k}\,\varphi^{\kern 0.4pt\alpha111}
=L_{\boldsymbol\Upsilon_k}(\varphi^{\kern 0.4pt\alpha111})
+\sum^2_{\theta=1}\uuA^\alpha_{k\kern 0.4pt\theta}
\,\varphi^{\kern 1.0pt\theta\kern 0.4pt 111}+3\,\uA^1_{k1}
\,\varphi^{\kern 0.6pt\alpha111}.
\mytag{2.2}
$$
Covariant differentiations for the wave functions of leptons and quarks
\mythetag{1.4}, \mythetag{1.5}, \mythetag{1.13}, and \mythetag{1.14} are
defined in a similar way:
$$
\allowdisplaybreaks
\gather
\hskip -2em
\gathered
\nabla_{\!k}\,\bulpsi^{a\alpha}_{111}[q]
=L_{\boldsymbol\Upsilon_k}(\bulpsi^{a\alpha}_{111}[q])
+\sum^4_{b=1}\Alpha^a_{kb}\,\bulpsi^{b\alpha}_{111}[q]\,+\\
+\sum^2_{\theta=1}\uuA^\alpha_{k\kern 0.4pt\theta}
\,\bulpsi^{a\kern 0.4pt\theta}_{111}[q]-3\,\uA^1_{k1}
\,\bulpsi^{a\alpha}_{111}[q]\text{, \ where \ }q=e,\mu,\tau;
\endgathered
\mytag{2.3}\\
\hskip -2em
\gathered
\nabla_{\!k}\,\circpsi^a_{111111}[q]
=L_{\boldsymbol\Upsilon_k}(\circpsi^a_{111111}[q])
+\sum^4_{b=1}\Alpha^a_{kb}\,\circpsi^b_{111111}[q]\,-\\
\vspace{1ex}
-\,6\,\uA^1_{k1}\,\circpsi^a_{111111}[q]\text{, \ where \ }q=e,\mu,\tau;
\endgathered
\mytag{2.4}\\
\hskip -2em
\gathered
\nabla_{\!k}\,\bulpsi^{a1\alpha\beta}[q]
=L_{\boldsymbol\Upsilon_k}(\bulpsi^{a1\alpha\beta}[q])
+\sum^4_{b=1}\Alpha^a_{kb}\,\bulpsi^{b1\alpha\beta}[q]
+\uA^1_{k1}\,\bulpsi^{a1\alpha\beta}[q]\,+\\
+\sum^2_{\theta=1}\uuA^\alpha_{k\kern 0.4pt\theta}\,
\bulpsi^{a1\theta\kern 0.4pt\beta}[q]+\sum^3_{\theta=1}
\uuuA^\alpha_{k\kern 0.4pt\theta}\,\bulpsi^{a1\alpha\kern 
0.4pt\theta}[q]\text{, \ where \ }q=1,2,3;
\endgathered
\mytag{2.5}\\
\hskip -2em
\gathered
\nabla_{\!k}\,\circpsi^{a1111\beta}[q]
=L_{\boldsymbol\Upsilon_k}(\circpsi^{a1111\beta}[q])
+\sum^4_{b=1}\Alpha^a_{kb}\,\circpsi^{b1111\beta}[q]\,+\\
+\,4\,\uA^1_{k1}\,\circpsi^{a1111\beta}[q]
+\sum^3_{\theta=1}\uuuA^\alpha_{k
\kern 0.4pt\theta}\,\circpsi^{a1111\kern 0.4pt\theta}[q]
\text{, \ where \ }q=u,c,t;
\endgathered
\mytag{2.6}\\
\hskip -2em
\gathered
\nabla_{\!k}\,\circpsi^{a\beta}_{11}[q]
=L_{\boldsymbol\Upsilon_k}(\circpsi^{a\beta}_{11}[q])
+\sum^4_{b=1}\Alpha^a_{kb}\,\circpsi^{b\beta}_{11}[q]\,-\\
-\,2\,\uA^1_{k1}\,\circpsi^{a\beta}_{11}[q]
+\sum^3_{\theta=1}\uuuA^\alpha_{k
\kern 0.4pt\theta}\,\circpsi^{a\kern 0.4pt\theta}_{11}[q]
\text{, \ where \ }q=d,s,b;
\endgathered
\mytag{2.7}
\endgather
$$
Here $\uuA^\alpha_{k\kern 0.4pt\theta}$ and $\uA^1_{k1}$ are the 
components of some connections associated with the electro-weak 
bundles $S\uuU\!M$ and $\uU\!M$, while $\uuuA^\alpha_{k\kern 0.4pt
\theta}$ are the components of some connection associated with the 
color bundle $S\uuuU\!M$. In Standard Model they are interpreted as 
the components of gauge fields.\par
     Apart from the gauge field connections, in \mythetag{2.3},
\mythetag{2.4}, \mythetag{2.5}, \mythetag{2.6}, and \mythetag{2.7} 
we see the components $\Alpha^a_{kb}$ of the spinor connection. These
quantities are due to the presence of a gravitation field. Unlike
$\uA^1_{k1}$, $\uuA^\alpha_{k\kern 0.4pt\theta}$, and $\uuuA^\alpha_{k
\kern 0.4pt\theta}$, the spinor connection components $\Alpha^a_{kb}$ 
are not interpreted as gauge fields. They are uniquely determined by 
the metric tensor $\bold g$ in the base space-time manifold $M$ (see
\mycite{7}).\par
    The gauge field connections $\uA$, $\uuA$, and $\uuuA$ are related 
to the basic fields in $\uU\!M$, $S\uuU\!M$, and $S\uuuU\!M$ through the
series of concordance conditions (see \mycite{6}):
$$
\xalignat 2
&\hskip -2em
\nabla\bolduD=0,
\mytag{2.8}
\\
&\hskip -2em
\nabla\bolduuD=0,
&&\nabla\bolduud=0,
\mytag{2.9}
\\
&\hskip -2em
\nabla\bolduuuD=0,
&&\nabla\bolduuud=0.
\mytag{2.10}
\endxalignat
$$
In a coordinate form the concordance conditions \mythetag{2.8},
\mythetag{2.9}, \mythetag{2.10} are written as
$$
\align
&\hskip -2em
\nabla_{\!k}\,\uD_{11}=
L_{\boldsymbol\Upsilon_k}(\uuD_{11})-\uD_{11}
\,\uA^1_{k1}-\uD_{11}
\,\overline{\uA^{\raise 0.6pt\hbox{$\ssize1$}}_{k1}}=0,
\mytag{2.11}\\
&\hskip -2em
\nabla_{\!k}\,\uuD_{i\bar j}=
L_{\boldsymbol\Upsilon_k}(\uuD_{i\bar j})-\sum^2_{a=1}\uuD_{a\bar j}
\,\uuA^a_{ki}-\sum^2_{\bar a=1}\uuD_{i\bar a}
\,\overline{\uuA^{\raise 0.6pt\hbox{$\ssize\bar a$}}_{k\bar j}}=0,
\mytag{2.12}\\
&\hskip -2em
\nabla_{\!k}\,\uud_{ij}=
L_{\boldsymbol\Upsilon_k}(\uud_{ij})-\sum^2_{a=1}
\uud_{aj}\,\uuA^a_{ki}-\sum^2_{a=1}\uud_{ia}\,\uuA^a_{kj}=0,
\mytag{2.13}\\
&\hskip -2em
\nabla_{\!k}\,\uuuD_{i\bar j}=
L_{\boldsymbol\Upsilon_k}(\uuuD_{i\bar j})-\sum^3_{a=1}\uuuD_{a\bar j}
\,\uuuA^a_{ki}-\sum^3_{\bar a=1}\uuuD_{i\bar a}
\,\overline{\uuuA^{\raise 0.6pt\hbox{$\ssize\bar a$}}_{k\bar j}}=0,
\mytag{2.14}\\
&\hskip -2em
\gathered
\nabla_{\!k}\,\uuud_{ijm}=
L_{\boldsymbol\Upsilon_k}(\uuud_{ijm})-\sum^3_{a=1}
\uuud_{ajm}\,\uuuA^a_{ki}-\sum^3_{a=1}\uuud_{iam}\,\uuuA^a_{kj}\,-\\
-\sum^3_{a=1}\uuud_{ija}\,\uuuA^a_{km}=0.
\endgathered
\mytag{2.15}
\endalign
$$
Apart from \mythetag{2.11}, \mythetag{2.12}, \mythetag{2.13},
\mythetag{2.14}, \mythetag{2.15}, the tensor fields $\bolduuD$ and
$\bolduuuD$ are related to $\bolduud$ and $\bolduuud$ by means of
the following concordance conditions (see \mycite{2}):
$$
\align
&\hskip -2em
\sum^2_{i=1}\sum^2_{j=1}\uud^{\kern 0.4pt ij}
\,\uuD_{i\bar i}\,\uuD_{j\bar j}
=-\overline{\uud_{\kern 0.4pt\bar i\kern 0.2pt\bar j}},
\mytag{2.16}\\
&\hskip -2em
\sum^3_{i=1}\sum^3_{j=1}\sum^3_{k=1}\uuud^{\kern 0.4pt ijk}
\,\uuuD_{i\bar i}\,\uuuD_{j\bar j}\,\uuuD_{k\bar k}
=\overline{\uuud_{\kern 0.4pt\bar i\kern 0.2pt\bar j\bar k}}.
\endalign
$$\par
     Each connection produces its curvature tensor. In the case of 
the gauge connections $\uA$, $\uuA$, and $\uuuA$ for the components
of the curvature tensors we have 
$$
\align
&\hskip -2em
\goth R^1_{1ij}=\goth R_{ij}=
L_{\boldsymbol\Upsilon_i}(\uA^1_{j1})-L_{\boldsymbol\Upsilon_j}(\uA^1_{i1})
-\sum^3_{h=0}c^{\,h}_{ij}\,\uA^1_{h1},
\mytag{2.17}\\
&\hskip -2em
\aligned
\goth R^p_{k\kern 0.4pt ij}&=L_{\boldsymbol\Upsilon_i}(\uuA^p_{j\,k})
-L_{\boldsymbol\Upsilon_j}(\uuA^p_{i\,k})\,+\\
&+\sum^2_{h=1}\left(\uuA^p_{i\,h}\,\uuA^{\!h}_{j\,k}
-\uuA^p_{j\,h}\,\uuA^{\!h}_{i\,k}\right)
-\sum^3_{h=0}c^{\,h}_{ij}\,\uuA^p_{hk},
\endaligned
\mytag{2.18}\\
&\hskip -2em
\aligned
\goth R^p_{k\kern 0.4pt ij}&=L_{\boldsymbol\Upsilon_i}(\uuuA^p_{j\,k})
-L_{\boldsymbol\Upsilon_j}(\uuuA^p_{i\,k})\,+\\
&+\sum^3_{h=1}\left(\uuuA^p_{i\,h}\,\uuuA^{\!h}_{j\,k}
-\uuuA^p_{j\,h}\,\uuuA^{\!h}_{i\,k}\right)
-\sum^3_{h=0}c^{\,h}_{ij}\,\uuuA^p_{hk}
\endaligned
\mytag{2.19}
\endalign
$$
(see \mycite{6}). Here $c^{\,h}_{ij}$ are the parameters introduced 
through the commutation relationships for the frame vector fields
$\boldsymbol\Upsilon_0$, $\boldsymbol\Upsilon_1$, $\boldsymbol
\Upsilon_2$, $\boldsymbol\Upsilon_3$:
$$
\hskip -2em
[\boldsymbol\Upsilon_{\!i},\,\boldsymbol\Upsilon_{\!j}]=
\sum^3_{k=0}c^{\,k}_{ij}\,\boldsymbol\Upsilon_{\!k}.
\mytag{2.20}
$$
From \mythetag{2.11}, \mythetag{2.12}, and \mythetag{2.14} we derive 
the following identities for the curvature tensors introduced by the 
formulas \mythetag{2.17}, \mythetag{2.18}, and \mythetag{2.19}:
$$
\align
&\hskip -2em
\goth R^1_{1ij}+\overline{\goth R^{\raise 0.6pt
\hbox{$\ssize 1$}}_{1ij}}=0,
\mytag{2.21}\\
&\hskip -2em
\sum^2_{\bar k=1}\goth R^k_{p\kern 0.6pt ij}\,\uuD_{k\bar q}
+\sum^2_{k=1}\overline{\goth R^{\raise 0.6pt
\hbox{$\ssize\bar k$}}_{\bar q\kern 0.2pt ij}}\,\uuD_{p\bar k}=0,
\mytag{2.22}\\
&\hskip -2em
\sum^3_{\bar k=1}\goth R^k_{p\kern 0.6pt ij}\,\uuuD_{k\bar q}
+\sum^3_{k=1}\overline{\goth R^{\raise 0.6pt
\hbox{$\ssize\bar k$}}_{\bar q\kern 0.2pt ij}}\,\uuuD_{p\bar k}=0.
\mytag{2.23}
\endalign
$$
Similarly, from \mythetag{2.13} and \mythetag{2.15} for 
\mythetag{2.18} and \mythetag{2.19} we derive
$$
\align
&\hskip -2em
\sum^2_{k=1}\goth R^p_{k\kern 0.2pt ij}\,\uud^{\kern 0.4pt kq}
+\sum^2_{k=1}\goth R^q_{k\kern 0.2pt ij}\,\uud^{\kern 0.4pt p
\kern 0.2pt k}=0,
\mytag{2.24}\\
&\hskip -2em
\sum^3_{k=1}\goth R^p_{k\kern 0.2pt ij}\,\uuud^{\kern 0.4pt kqm}
+\sum^3_{k=1}\goth R^q_{k\kern 0.2pt ij}\,\uuud^{\kern 0.4pt p
\kern 0.2pt km}
+\sum^3_{k=1}\goth R^m_{k\kern 0.2pt ij}\,\uuud^{\kern 0.4pt p
\kern 0.2pt qk}=0.
\mytag{2.25}
\endalign
$$
The equality \mythetag{2.21} means that $\goth R^1_{1ij}$ given
by the formula \mythetag{2.17} is a purely imaginary number. The
equalities \mythetag{2.22} and \mythetag{2.23} are more complicated.
They mean that for fixed $i$ and $j$ the components of the curvature
tensors \mythetag{2.18} and \mythetag{2.19} form skew-Hermitian 
matrices with respect to Hermitian forms $\bolduuD$ and $\bolduuuD$
respectively.\par
    The equality \mythetag{2.25} looks rather complicated. However, 
for each fixed $i$ and $j$ it is equivalent to the following zero 
trace condition for the curvature tensor \mythetag{2.19}:
$$
\hskip -2em
\sum^3_{k=1}\goth R^k_{k\kern 0.2pt ij}=0.
\mytag{2.26}
$$
From \mythetag{2.24} one can derive the analogous equality for
\mythetag{2.18}:
$$
\hskip -2em
\sum^2_{k=1}\goth R^k_{k\kern 0.2pt ij}=0.
\mytag{2.27}
$$
Unlike \mythetag{2.26}, in this case \mythetag{2.24} is not equivalent 
to \mythetag{2.27}. Remember that in $S\uuU\!M$ the tensor $\bolduud$ 
is used for raising and lowering indices. Let's denote
$$
\xalignat 2
&\hskip -2em
\goth R^{p\kern 0.6pt q}_{ij}
=\sum^2_{k=1}\goth R^p_{k\kern 0.2pt ij}\,\uud^{\kern 0.4pt kq},
&&\goth R_{p\kern 0.6pt q\kern 0.2pt ij}
=\sum^2_{k=1}\goth R^k_{q\kern 0.2pt ij}\,\uud_{kq}.
\mytag{2.28}
\endxalignat
$$
In terms of \mythetag{2.28} the equality \mythetag{2.24} is equivalent
to the symmetry conditions
$$
\xalignat 2
&\goth R^{p\kern 0.6pt q}_{ij}=\goth R^{q\kern 0.4pt p}_{ij},
&&\goth R_{p\kern 0.6pt q\kern 0.2pt ij}=\goth R_{qp\kern 0.2pt ij}.
\endxalignat
$$\par
     The curvature tensors \mythetag{2.17}, \mythetag{2.18}, 
\mythetag{2.19} are physical fields associated with the gauge fields
$\uA$, $\uuA$ and $\uuuA$. There are the following identities
for them:
$$
\align
&\hskip -2em
\uD^{11}\,\uD_{11}\,\goth R^1_{1ij}
\,\overline{\goth R^{\raise 0.6pt\hbox{$\ssize1$}}_{1mn}}=
-\goth R^1_{1ij}\,\goth R^1_{1mn},
\mytag{2.29}\\
\vspace{1.5ex}
&\hskip -2em
\sum^2_{p=1}\sum^2_{\bar p=1}\sum^2_{q=1}\sum^2_{\bar q=1}
\uuD^{q\bar q}\,\uuD_{p\bar p}\ \goth R^p_{q\kern 0.2pt ij}
\ \overline{\goth R^{\raise 0.6pt
\hbox{$\ssize\bar p$}}_{\bar q\kern 0.2pt mn}}
=-\sum^2_{p=1}\sum^2_{q=1}
\goth R^p_{q\kern 0.2pt ij}\ \goth R^q_{p\kern 0.2pt ij},
\mytag{2.30}\\
&\hskip -2em
\sum^3_{p=1}\sum^3_{\bar p=1}\sum^3_{q=1}\sum^3_{\bar q=1}
\uuuD^{q\bar q}\,\uuuD_{p\bar p}\ \goth R^p_{q\kern 0.2pt ij}
\ \overline{\goth R^{\raise 0.6pt
\hbox{$\ssize\bar p$}}_{\bar q\kern 0.2pt mn}}
=-\sum^3_{p=1}\sum^3_{q=1}
\goth R^p_{q\kern 0.2pt ij}\ \goth R^q_{p\kern 0.2pt ij}.
\mytag{2.31}
\endalign
$$
By $\uD^{11}$ in \mythetag{2.29} we denote the component of the inverse
Hermitian metric tensor of the bundle $\uU\!M$. It is related to $\uD_{11}$
by means of the formula
$$
\hskip -2em
\uD^{11}=\frac{1}{\uD_{11}}.
\mytag{2.32}
$$
Similarly, by $\uuD^{q\bar q}$ and $\uuuD^{q\bar q}$ in \mythetag{2.30}
and \mythetag{2.31} we denote the components of the inverse Hermitian 
metric tensors in $S\uuU\!M$ and $S\uuuU\!M$. They form two matrices
inverse to the matrices $\uuD_{p\bar p}$ and $\uuuD_{p\bar p}$ respectively:
$$
\xalignat 2
&\hskip -2em
\sum^2_{\bar p=1}\uuD_{p\bar p}\,\uuD^{q\bar p}=\delta^q_p,
&&\sum^2_{p=1}\uuD_{p\bar p}\,\uuD^{p\bar q}=\delta^{\bar q}_{\bar p},
\mytag{2.33}\\
&\hskip -2em
\sum^3_{\bar p=1}\uuuD_{p\bar p}\,\uuuD^{q\bar p}=\delta^q_p,
&&\sum^3_{p=1}\uuuD_{p\bar p}\,\uuuD^{p\bar q}=\delta^{\bar q}_{\bar p},
\mytag{2.34}
\endxalignat
$$
The formula \mythetag{2.29} follows immediately from 
\mythetag{2.21} and \mythetag{2.32}. The formula \mythetag{2.30}
is derived from \mythetag{2.22} with the use of the formulas
\mythetag{2.33}. Similarly, the formula \mythetag{2.31}
is derived from \mythetag{2.23} with the use of the formulas
\mythetag{2.34}.
\mydefinition{2.1} The electro-weak and color bundles $\uU\!M$,
$S\uuU\!M$, and $S\uuuU\!M$ are called {\it flat\/} if there 
are three flat connections $\uA$, $\uuA$, and $\uuuA$ in these
bundles, i\.\,e\. three connections with zero curvature tensors
\pagebreak \mythetag{2.17}, \mythetag{2.18}, and \mythetag{2.19}.
\enddefinition
     The connections $\uA$, $\uuA$, and $\uuuA$ in the above
definition~\mythedefinition{2.1} are assumed to be concordant
with the basic fields $\bolduD$, $\bolduuD$, $\bolduuuD$, 
$\bolduud$, and $\bolduuud$ of the electro-weak and color 
bundles. Therefore, the flatness condition
$$
\hskip -2em
\eufb R=0
\mytag{2.35}
$$
written for their curvature tensors \mythetag{2.17}, \mythetag{2.18}, 
and \mythetag{2.19} is an additional condition complementary to
\mythetag{2.8}, \mythetag{2.9}, and \mythetag{2.10}. Below we implicitly
assume that the bundles $\uU\!M$, $S\uuU\!M$, and $S\uuuU\!M$ are
flat in the sense of the definition~\mythedefinition{2.1}.
\mydefinition{2.2} The components of the flat connections $\uA$, $\uuA$,
and $\uuuA$ constitute a classical vacuum of bosonic gauge fields of the
Standard Model.
\enddefinition
The existence of a classical gauge vacuum is provided by the 
definition~\mythedefinition{2.1}. The problem of uniqueness of such
a vacuum in the case of flat bundles as well as the possibility of
non-flat vacua should be studied in separate papers.\par
     From now on let's assume that some flat electro-weak and color
bundles $\uU\!M$, $S\uuU\!M$, and $S\uuuU\!M$ are chosen and some flat
gauge vacuum of them is fixed. Then for general non-vacuum gauge fields 
we write
$$
\align
&\hskip -2em
\uA^1_{k1}=\uA^1_{k1}[vac]-
\frac{i\,e\,g_1}{\hbar\,c}\,\uCA^1_{k1},
\mytag{2.36}
\\
&\hskip -2em
\uuA^\alpha_{k\beta}=\uuA^\alpha_{k\beta}[vac]-
\frac{i\,e\,g_2}{\hbar\,c}\,\uuCA^\alpha_{k\beta}.
\mytag{2.37}\\
&\hskip -2em
\uuuA^\alpha_{k\beta}=\uuuA^\alpha_{k\beta}[vac]-
\frac{i\,e\,g_3}{\hbar\,c}\,\uuuCA^\alpha_{k\beta}.
\mytag{2.38}
\endalign
$$
Here $e$ is the charge of an electron\footnote{\ Note that $e$ in
\mythetag{2.36}, \mythetag{2.37}, and \mythetag{2.38} is a positive
quantity, therefore, it is better to say that $e$ is the charge of
a positron. The numeric values \mythetag{2.39} of the physical constants
$e$, $\hbar$, and $c$ are taken from the NIST site
\myhref{http://physics.nist.gov/cuu/Constants}
{http:/\negskp/physics.nist.gov/cuu/Constants}. The value of $e$ there
is given in SI units. It is converted to CGS units by means of the
NIST value of $\varepsilon_0$.}, $c$ is the light speed, and $\hbar$ 
is the Planck constant. Below are the numeric values of these 
foundamental constants:
\adjustfootnotemark{-1}
$$
\align
&\hskip -2em
e\approx 4.80420440\cdot 10^{-10}\,\text{\it g}^{\,\sssize 1/2}
\cdot\text{\it cm}^{\,\sssize 3/2}\cdot\text{\it sec}^{\sssize\,-1}
,\\
&\hskip -2em
\hbar\approx 1.05457168\cdot 10^{-27}\,\text{\it g}\cdot
\text{\it cm}^{\,\sssize 2}\cdot\text{\it sec}^{\sssize\,-1},
\mytag{2.39}\\
&\hskip -2em
c\approx 2.99792458\cdot 10^{10}\,\text{\it cm}\cdot
\text{\it sec}^{\sssize\,-1}.
\endalign
$$
By $g_1$, $g_2$, and $g_3$ in \mythetag{2.36}, \mythetag{2.37}, and
\mythetag{2.38} we denote some numeric constants which are called the
{\it coupling constants}. Their values are determined a posteriori by
comparing the Standard Model predictions with experimental data.\par
     Unlike $\uA^\alpha_{k\beta}$, $\uuA^\alpha_{k\beta}$, and 
$\uuuA^\alpha_{k\beta}$, the quantities $\uCA^\alpha_{k\beta}$,
$\uuCA^\alpha_{k\beta}$, and $\uuuCA^\alpha_{k\beta}$ in \mythetag{2.36},
\mythetag{2.37}, and \mythetag{2.38} are the components of tensor fields.
We denote these tensor fields by $\bolduCA$, $\bolduuCA$, and $\bolduuuCA$
respectively. Since $\eufb R[vac]=0$, substituting \mythetag{2.36}, 
\mythetag{2.37}, and \mythetag{2.38} into \mythetag{2.17},
\mythetag{2.18}, and \mythetag{2.19} respectively, we derive
$$
\allowdisplaybreaks
\align
&\hskip -6em
\goth R^1_{1\kern 0.2pt ij}=-\frac{i\,e\,g_1}{\hbar\,c}
\left(\nabla_{\!i}\,\uCA^{1\vphantom{p}}_{j\,1}-\nabla_{\!j}
\,\uCA^1_{i\,1}\right)
\mytag{2.40}\\
&\hskip -6em
\goth R^p_{k\kern 0.2pt ij}=-\frac{i\,e\,g_2}{\hbar\,c}
\left(\nabla_{\!i}\,\uuCA^p_{j\,k}-\nabla_{\!j}\,\uuCA^p_{i\,k}\right)
-\sum^2_{h=1}\left(\frac{e\,g_2}{\hbar\,c}\right)^{\!2}\!
\left(\uuCA^p_{i\,h}\,\uuCA^{\!h}_{j\,k}
-\uuCA^p_{j\,h}\,\uuCA^{\!h}_{i\,k}\right),
\hskip -2em
\mytag{2.41}\\
&\hskip -6em
\goth R^p_{k\kern 0.2pt ij}=-\frac{i\,e\,g_3}{\hbar\,c}
\left(\nabla_{\!i}\,\uuCA^p_{j\,k}-\nabla_{\!j}\,\uuCA^p_{i\,k}\right)
-\sum^3_{h=1}\left(\frac{e\,g_3}{\hbar\,c}\right)^{\!2}\!
\left(\uuCA^p_{i\,h}\,\uuCA^{\!h}_{j\,k}
-\uuCA^p_{j\,h}\,\uuCA^{\!h}_{i\,k}\right).
\hskip -2em
\mytag{2.42}
\endalign
$$
Here $\nabla_{\!i}$ and $\nabla_{\!j}$ are vacuum covariant derivatives,
i\.\,e\. they are calculated with respect to the vacuum electro-weak
and color connections. In deriving \mythetag{2.40}, \mythetag{2.41},
and \mythetag{2.42} from \mythetag{2.17}, \mythetag{2.18}, and
\mythetag{2.19} we used the following equality for the components of the
Levi-Civita connection:
$$
\hskip -2em
\Gamma^h_{ij}-\Gamma^h_{j\kern 0.4pt i}=c^{\,h}_{ij}.
\mytag{2.43}
$$
The equality \mythetag{2.43} means that $\Gamma^h_{ij}$ are the 
components of a torsion-free (symmetric) connection (see \thetag{2.8} 
in \mycite{6}). The quantities $c^{\,h}_{ij}$ are introduced by the
formula \mythetag{2.20}.\par
     Relying upon the formulas \mythetag{2.40}, \mythetag{2.41}, and
\mythetag{2.42}, now we introduce three tensor fields $\bolduF$, $\bolduuF$,
and $\bolduuuF$ with the following components:
$$
\align
&\hskip -2em
\uF_{\!\!ij}=\nabla_{\!i}\,\uCA^1_{j\,1}-\nabla_{\!j}\,
\uCA^1_{i\,1},
\mytag{2.44}\\
&\hskip -2em
\uuF^p_{\!\!kij}=\nabla_{\!i}\,\uuCA^p_{j\,k}-\nabla_{\!j}\,
\uuCA^p_{i\,k}-\frac{i\,e\,g_2}{\hbar\,c}\sum^2_{h=1}\left(\uuCA^p_{i\,h}
\,\uuCA^{\!h}_{j\,k}-\uuCA^p_{j\,h}\,\uuCA^{\!h}_{i\,k}\right),
\mytag{2.45}\\
&\hskip -2em
\uuuF^p_{\!\!kij}=\nabla_{\!i}\,\uuuCA^p_{j\,k}-\nabla_{\!j}\,
\uuuCA^p_{i\,k}-\frac{i\,e\,g_3}{\hbar\,c}\sum^3_{h=1}\left(\uuuCA^p_{i\,h}
\,\uuuCA^{\!h}_{j\,k}-\uuuCA^p_{j\,h}\,\uuuCA^{\!h}_{i\,k}\right).
\mytag{2.46}
\endalign
$$
Due to \mythetag{2.44}, \mythetag{2.45}, and \mythetag{2.46} the formulas
\mythetag{2.40}, \mythetag{2.41}, \mythetag{2.42} are written as
$$
\xalignat 3
&\hskip -1em
\goth R^1_{1\kern 0.2pt ij}=-\frac{i\,e\,g_1}{\hbar\,c}
\,\uF_{\!\!ij},
&&\goth R^p_{k\kern 0.2pt ij}=-\frac{i\,e\,g_2}{\hbar\,c}
\,\uuF^p_{\!\!kij},
&&\goth R^p_{k\kern 0.2pt ij}=-\frac{i\,e\,g_3}{\hbar\,c}
\,\uuuF^p_{\!\!kij},\qquad\quad
\mytag{2.47}
\endxalignat
$$\par
     The concordance conditions \mythetag{2.8}, \mythetag{2.9}, and
\mythetag{2.10} should be fulfilled for both vacuum and non-vacuum 
connections in \mythetag{2.36}, \mythetag{2.37}, and \mythetag{2.38}.
Therefore, from the formulas \mythetag{2.11}, \mythetag{2.12}, and
\mythetag{2.13} we derive
$$
\align
&\hskip -2em
\uCA^1_{k1}=\overline{\uCA^{\raise 0.6pt\hbox{$\ssize 1$}}_{k1}},
\mytag{2.48}\\
&\hskip -2em
\sum^2_{a=1}\uuD_{a\bar j}
\,\uuCA^a_{ki}=\sum^2_{\bar a=1}\uuD_{i\bar a}
\,\overline{\uuCA^{\kern 0.2pt\raise 0.6pt
\hbox{$\ssize\bar a$}}_{k\bar j}},
\mytag{2.49}\\
&\hskip -2em
\sum^2_{a=1}\uuCA^a_{ki}\,\uud_{aj}
=\sum^2_{a=1}\uuCA^a_{kj}\,\uud_{ai}.
\mytag{2.50}
\endalign
$$
Similarly, from the formulas \mythetag{2.14} and \mythetag{2.15} we derive
$$
\align
&\hskip -2em
\sum^3_{a=1}\uuuD_{a\bar j}
\,\uuuCA^a_{ki}=\sum^3_{\bar a=1}\uuuD_{i\bar a}
\,\overline{\uuuCA^{\kern 0.2pt\raise 0.6pt
\hbox{$\ssize\bar a$}}_{k\bar j}},
\mytag{2.51}\\
&\hskip -2em
\sum^3_{a=1}
\uuud_{ajm}\,\uuuCA^a_{ki}+\sum^3_{a=1}\uuud_{iam}\,\uuuCA^a_{kj}
+\sum^3_{a=1}\uuud_{ija}\,\uuuCA^a_{km}=0.
\mytag{2.52}
\endalign
$$
Now let's substitute \mythetag{2.47} into \mythetag{2.21}, 
\mythetag{2.22}, \mythetag{2.23}, \mythetag{2.24}, and \mythetag{2.25}.
\pagebreak As a result we get a series of formulas similar to 
\mythetag{2.48}, \mythetag{2.49}, \mythetag{2.50}, \mythetag{2.51}, 
and \mythetag{2.52}:
$$
\align
&\hskip -2em
\uF_{\!\!ij}=\overline{\uF_{\!\!ij}},
\mytag{2.53}\\
&\hskip -2em
\sum^2_{\bar k=1}\uuF^k_{\!\!p\kern 0.6pt ij}\,\uuD_{k\bar q}
=\sum^2_{k=1}\overline{\uuF^{\raise 0.6pt \hbox{$\ssize\bar k$}}_{\!
\!\bar q\kern 0.2pt ij}}\,\uuD_{p\bar k},
\mytag{2.54}\\
&\hskip -2em
\sum^3_{\bar k=1}\uuuF^k_{\!\!p\kern 0.6pt ij}\,\uuuD_{k\bar q}
=\sum^3_{k=1}\overline{\uuuF^{\raise 0.6pt
\hbox{$\ssize\bar k$}}_{\!\!\bar q\kern 0.2pt ij}}\,\uuuD_{p\bar k}.
\mytag{2.55}\\
&\hskip -2em
\sum^2_{k=1}\uuF^k_{\!\!p\kern 0.6pt ij}\,\uud_{kq}
=\sum^2_{k=1}\uuF^k_{\!\!q\kern 0.6pt ij}\,\uud_{kp},\\
&\hskip -2em
\sum^3_{k=1}\uuud_{kqm}\,\uuuF^k_{\!\!p\kern 0.6pt ij}
+\sum^3_{k=1}\uuud_{pkm}\,\uuuF^k_{\!\!q\kern 0.6pt ij}
+\sum^3_{k=1}\uuud_{p\kern 0.4pt qk}\,\uuuF^k_{\!\!q\kern 0.6pt ij}=0.
\endalign
$$
Substituting \mythetag{2.47} into \mythetag{2.29}, \mythetag{2.30},
and \mythetag{2.31}, we derive the following identities:
$$
\align
&\hskip -2em
\uD^{11}\,\uD_{11}\ \uF_{\!\!ij}
\ \overline{\uF_{\!\!mn}}
=\uF_{\!\!ij}\ \uF_{\!\!mn},
\mytag{2.56}\\
\vspace{1.5ex}
&\hskip -2em
\sum^2_{p=1}\sum^2_{\bar p=1}\sum^2_{q=1}\sum^2_{\bar q=1}
\uuD^{q\bar q}\,\uuD_{p\bar p}\ \uuF^p_{\!\!q\kern 0.2pt ij}
\ \overline{\uuF^{\raise 0.6pt
\hbox{$\ssize\bar p$}}_{\!\!\bar q\kern 0.2pt mn}}
=\sum^2_{p=1}\sum^2_{q=1}\uuF^p_{\!\!q\kern 0.2pt ij}
\ \uuF^q_{\!\!p\kern 0.2pt ij},
\mytag{2.57}\\
&\hskip -2em
\sum^3_{p=1}\sum^3_{\bar p=1}\sum^3_{q=1}\sum^3_{\bar q=1}
\uuuD^{q\bar q}\,\uuuD_{p\bar p}\ \uuuF^p_{\!\!q\kern 0.2pt ij}
\ \overline{\uuuF^{\raise 0.6pt
\hbox{$\ssize\bar p$}}_{\!\!\bar q\kern 0.2pt mn}}
=\sum^3_{p=1}\sum^3_{q=1}\uuuF^p_{\!\!q\kern 0.2pt ij}
\ \uuuF^q_{\!\!p\kern 0.2pt ij}.
\mytag{2.58}
\endalign
$$
The identities \mythetag{2.56}, \mythetag{2.57}, \mythetag{2.58} 
can be derived directly from  \mythetag{2.53}, \mythetag{2.54},
\mythetag{2.55} with the use of the formulas \mythetag{2.32},
\mythetag{2.33}, and \mythetag{2.34}.\par
    Let's denote by $\Cal L_1=\Cal L_1(\bolduF)$, $\Cal L_2
=\Cal L_2(\bolduuF)$, and $\Cal L_3=\Cal L_3(\bolduuuF)$ the kinetic 
terms of the gauge fields in the action integral of the Standard Model:
$$
\align
&\hskip -2em
\Cal L_1=-\frac{1}{16\,\pi\,c}
\int\msum{3}\Sb i,j,m,n=0\endSb g^{im}\,g^{jn}\,
\uF^1_{\!\!1ij}\ \uF^1_{\!\!1mn}\,dV,
\mytag{2.59}\\
&\hskip -2em
\Cal L_2=-\frac{1}{32\,\pi\,c}
\int\sum^2_{p=1}\sum^2_{q=1}
\msum{3}\Sb i,j,m,n=0\endSb g^{im}\,g^{jn}\,
\uuF^p_{\!\!q\kern 0.2pt ij}
\ \uuF^q_{\!\!p\kern 0.2pt ij}\,dV,
\mytag{2.60}\\
&\hskip -2em
\Cal L_3=-\frac{1}{48\,\pi\,c}
\int\sum^3_{p=1}\sum^3_{q=1}
\msum{3}\Sb i,j,m,n=0\endSb g^{im}\,g^{jn}\,
\uuuF^p_{\!\!q\kern 0.2pt ij}
\ \uuuF^q_{\!\!p\kern 0.2pt ij}\,dV.
\mytag{2.61}
\endalign
$$
Here $dV$ is the $4$-dimensional volume element in the base space-time
manifold $M$:
$$
\hskip -2em
dV=\sqrt{-\det\bold g\,}\,d^{\kern 0.5pt 4}\kern -0.5pt x=
\sqrt{-\det\bold g\,}\ dx^0\!\wedge dx^1\!\wedge dx^2\!\wedge dx^3.
\mytag{2.62}
$$
Due to the identities \mythetag{2.56}, \mythetag{2.57}, \mythetag{2.58} 
the integrals \mythetag{2.59}, \mythetag{2.60}, and \mythetag{2.61} are
real numbers. The coefficients preceding these integrals in \mythetag{2.59},
\mythetag{2.60}, and \mythetag{2.61} are chosen by analogy to the
Electrodynamics (see \mycite{8}).
\mydefinition{2.3} A Higgs field \mythetag{2.1} is called a flat classical
Higgs vacuum if 
$$
\xalignat 2
&\hskip -2em
\sum^2_{\alpha=1}\sum^2_{\bar\alpha=1}\uuD_{\alpha\bar\alpha}\,
\uD_{11}\,\uD_{11}\,\uD_{11}\,\varphi^{\kern 0.4pt\alpha111}\,
\overline{\varphi^{\kern 0.4pt\raise 
1.2pt\hbox{$\ssize\bar\alpha111$}}}=\const,
&&\nabla_{\!k}\varphi^{\kern 0.4pt\alpha111}=0,\qquad\quad
\mytag{2.63}
\endxalignat
$$
where the covariant derivatives $\nabla_{\!k}$ (see \mythetag{2.2}) are
determined by the flat connections forming a classical vacuum of gauge
fields in the sense of the definition~\mythedefinition{2.2}.
\enddefinition
    Due to \mythetag{2.8} and \mythetag{2.9} the conditions \mythetag{2.63}
are consistent with each other. Note that the Higgs field $\varphi^{\kern
0.4pt\alpha111}$ has neither spinor nor tensorial indices associated with
the tangent bundle $TM$. It is a scalar field with respect to Lorentz
transformations. For this reason the flatness condition \mythetag{2.35} 
fulfilled for the vacuum gauge fields guarantees the existence of local
nonzero Higgs vacua. As for a global nonzero Higgs vacuum, its existence
depends on the topology of the base manifold $M$ and the electro-weak 
bundles $\uU\!M$ and $S\uuU\!M$ over this base. Below we implicitly assume
that at least one nonzero classical Higgs vacuum does exist. We denote 
it by $\boldsymbol\varphi[vac]$.\par
    For any given Higgs field \mythetag{2.1} we denote by 
$|\boldsymbol\varphi|^2$ the left hand side of the first formula
\mythetag{2.63} in the above definition~\mythedefinition{2.3}:
$$
\hskip -2em
|\boldsymbol\varphi|^2=
\sum^2_{\alpha=1}\sum^2_{\bar\alpha=1}\uuD_{\alpha\bar\alpha}\,
\uD_{11}\,\uD_{11}\,\uD_{11}\,\varphi^{\kern 0.4pt\alpha111}\,
\overline{\varphi^{\kern 0.4pt\raise 
1.2pt\hbox{$\ssize\bar\alpha111$}}}.
\mytag{2.64}
$$
In terms of \mythetag{2.64} the potential of the Higgs field is written as
$$
\hskip -2em
V(\boldsymbol\varphi)=\lambda\,|\boldsymbol\varphi|^4
-\mu^2\,|\boldsymbol\varphi|^2
\mytag{2.65}
$$
(see \mycite{3}). Here $\lambda$ and $\mu$ are two constants. They 
are parameters of the Standard Model. By analogy to \mythetag{2.64} 
we denote by $|\nabla\boldsymbol\varphi|^2$ the following expression:
$$
\hskip -2em
|\nabla\boldsymbol\varphi|^2
=\sum^3_{i=0}\sum^3_{j=0}\sum^2_{\alpha=1}\sum^2_{\bar\alpha=1}
g^{ij}\,\uuD_{\alpha\bar\alpha}\,\uD_{11}\,\uD_{11}\,\uD_{11}\,
\nabla_{\!i}\,\varphi^{\kern 0.4pt\alpha 111}\,
\overline{\nabla_{\!j}\,\varphi^{\kern 0.4pt\raise 1.2pt
\hbox{$\ssize\bar\alpha 111$}}}.
\mytag{2.66}
$$
The expressions \mythetag{2.65} and  \mythetag{2.66} determine the
potential and kinetic terms for the Higgs field in the action integral
of the Standard Model:
$$
\align
&\hskip -2em
\Cal L_4=\frac{\hbar^2}{2\,m_\varphi\,c}\int|\nabla\boldsymbol\varphi|^2\,dV,
\mytag{2.67}\\
\vspace{2ex}
&\hskip -2em
\Cal L_5=-\frac{m_\varphi\,c}{2}\int V(\boldsymbol\varphi)\,dV.
\mytag{2.68}
\endalign
$$
Here $|\nabla\boldsymbol\varphi|^2$ and $V(\boldsymbol\varphi)$ are
given by the formulas \mythetag{2.66} and \mythetag{2.65}, while
$dV$ is the $4$-dimensional volume element given by the formula
\mythetag{2.62}. Like $\lambda$ amd $\mu$, the constant $m_\varphi$ 
is a parameter of the Standard Model, it is interpreted as the mass
of the Higgs field. Let's take the restricted action 
$$
\pagebreak
\hskip -2em
\Cal L_\varphi=\Cal L_4+\Cal L_5.
\mytag{2.69}
$$
The extremum of the restricted action \mythetag{2.69} is determined 
by the following Klein-Gordon-Fock type equation for the scalar Higgs
field $\boldsymbol\varphi$:
$$
-\frac{\hbar^2}{2\,m_\chi\,c}\sum^3_{i=0}\sum^3_{j=0}g^{ij}
\,\nabla_{\!i}\nabla_{\!j}\,\varphi^{\kern 0.4pt\alpha111}
+\frac{m_\chi\,c}{2}\left(2\,\lambda\,|\boldsymbol\varphi|^2
\,-\mu^2\right)\varphi^{\kern 0.4pt\alpha111}=0.
\qquad
\mytag{2.70}
$$
The equation \mythetag{2.70} has the trivial symmetric solution
$$
\hskip -2em
\boldsymbol\varphi=0.
\mytag{2.71}
$$
However, apart from this trivial solution \mythetag{2.71}, we shall
consider a nontrivial global solution described by the above
definition~\mythedefinition{2.3}:
$$
\hskip -2em
\boldsymbol\varphi=\boldsymbol\varphi[vac].
\mytag{2.72}
$$
Its existence is an assumption of the Standard Model when it is implemented
on a non-flat space-time manifold $M$ in General Relativity. Applying 
\mythetag{2.63} to the equation \mythetag{2.70}, for the nontrivial Higgs
vacuum we find
$$
\hskip -2em
|\boldsymbol\varphi[vac]|^2=\frac{\mu^2}{2\,\lambda}=\const.
\mytag{2.73}
$$
The constant \mythetag{2.73} is denoted be $v^2/2$. Then we can express
$\lambda$ through $\mu$ and $v$:
$$
\xalignat 2
&\hskip -2em
|\boldsymbol\varphi[vac]|^2=\frac{v^2}{2},
&&\lambda=\frac{\mu^2}{v^2}.
\qquad
\mytag{2.74}
\endxalignat
$$
\par
    Having defined the gauge and Higgs vacuum fields, now we proceed to
other matter fields, i\.\,e\. to leptons and quarks. Vacuum values of
their fields are trivial:
$$
\xalignat 3
&\hskip -2em
\boldsymbol\psi[e][vac]=0,
&&\boldsymbol\psi[\mu][vac]=0,
&&\boldsymbol\psi[\tau][vac]=0,
\qquad
\mytag{2.75}\\
&\hskip -2em
\boldsymbol\psi[1][vac]=0,
&&\boldsymbol\psi[2][vac]=0,
&&\boldsymbol\psi[3][vac]=0
\qquad
\mytag{2.76}
\endxalignat
$$
(compare \mythetag{2.75} and \mythetag{2.76} with \mythetag{1.3} and
\mythetag{1.12}). The formulas \mythetag{2.75} and \mythetag{2.76} mean
that both chiral and antichiral components of vacuum fields \mythetag{1.4},
\mythetag{1.5}, \mythetag{1.13}, and \mythetag{1.14} are zero.\par
\head
Gauge transformations\\
and perturbations of the Higgs vacuum.
\endhead
     Unlike the trivial Higgs vacuum \mythetag{2.71}, a nontrivial vacuum
\mythetag{2.72} is not unique. Applying some gauge transformation to a 
given one, we can get another nontrivial Higgs vacuum equivalent to it.
Assume that some orthonormal frame $(U,\,\bolduuPsi_1,\,\bolduuPsi_2)$ 
of the bundle $S\uuU\!M$ is fixed (see definition in \mycite{2}). Then the
basic fields of this bundle $\bolduuD$ and $\bolduud$ are given by the 
following constant matrices:
$$
\xalignat 2
&\hskip -2em
\uuD_{i\bar j}=\Vmatrix 1 & 0\\0 & 1\endVmatrix,
&&\uud_{ij}=\Vmatrix 0 & 1\\-1 & 0\endVmatrix.
\mytag{3.1}
\endxalignat
$$
The concordance conditions \mythetag{2.12} and \mythetag{2.13} written
for \mythetag{3.1} then mean that for each fixed $k$ the connection
components $\uuA^a_{ki}$ form a matrix from the Lie algebra $\MatAlgSU(2)$
of the special unitary group $\MatGrSU(2)$. This fact reflects the 
$\MatGrSU(2)$ symmetry of the Standard Model. Now assume that 
$\Omega^{\,\alpha}_\beta(p)$ is a matrix-valued smooth function with
the argument $p\in M$ and with values $\Omega^{\,\alpha}_\beta(p)\in
\MatGrSU(2)$. Let's denote $[\Omega^{-1}(p)]^{\,\alpha}_\beta$ the
components of the inverse matrix for $\Omega^{\,\alpha}_\beta(p)$. Then 
for each fixed $k$ the following quantities form a matrix from the Lie 
algebra $\MatAlgSU(2)$:
$$
\hskip -2em
\omega^{\,a}_{ki}=\sum^2_{s=1}
L_{\boldsymbol\Upsilon_k}(\Omega^a_s(p))
\,[\Omega^{-1}(p)]^s_i.
\mytag{3.2}
$$
The $\MatGrSU(2)$ gauge transformation for the Higgs field $\boldsymbol
\varphi$ is given by the formulas
$$
\align
\hskip -2em
\varphi^{\kern 0.4pt\alpha 111}&\longmapsto\sum^2_{\beta=1}
\Omega^{\,\alpha}_\beta\,\varphi^{\kern 0.4pt\beta 111}\!,
\mytag{3.3}\\
\hskip -2em
\uuA^a_{ki}&\longmapsto\sum^2_{s=1}\sum^2_{b=1}\Omega^a_b(p)
\,\uuA^b_{ks}\,[\Omega^{-1}(p)]^s_i-\omega^{\,a}_{ki},
\mytag{3.4}
\endalign
$$
where $\omega^{\,a}_{ki}$ is expressed through $\Omega^{\,\alpha}_\beta$
according to the formula \mythetag{3.2}. Using \mythetag{3.3} and
\mythetag{3.4} and taking into account \mythetag{3.2}, we derive
$$
\align
\hskip -2em
\nabla_{\!k}\varphi^{\kern 0.4pt\alpha 111}&\longmapsto\sum^2_{\beta=1}
\Omega^{\,\alpha}_\beta\,\nabla_{\!k}\varphi^{\kern 0.4pt\beta 111},
\mytag{3.5}\\
\nabla_{\!i}\nabla_{\!j}\,\varphi^{\kern 0.4pt\alpha 111}&\longmapsto
\sum^2_{\beta=1}\Omega^{\,\alpha}_\beta\,\nabla_{\!i}\nabla_{\!j}
\,\varphi^{\kern 0.4pt\alpha111}.
\mytag{3.6}
\endalign
$$
Due to \mythetag{3.5} and \mythetag{3.6}, applying a gauge transformation
\mythetag{3.3} to some solution of the equation \mythetag{2.70}, we get
a solution of the similar equation with respect to the transformed
connection \mythetag{3.4}. It is important to note that the initial and
transformed connections in \mythetag{3.4} have the equivalent curvature
tensors, i\.\,e\. we can complement the formulas \mythetag{3.3} and
\mythetag{3.4} with the formula
$$
\goth R^q_{k\kern 0.4pt ij}\longmapsto
\sum^2_{r=1}\sum^2_{s=1}\Omega^q_r(p)
\,\goth R^r_{s\kern 0.4pt ij}\,[\Omega^{-1}(p)]^s_k.
\mytag{3.7}
$$
Due to the formula \mythetag{3.7}, if the initial gauge vacuum in
\mythetag{2.37} is flat, then the transformed gauge vacuum obtained
by means of the formula \mythetag{3.4} is also flat.\par
     Now let's consider a perturbation of the Higgs vacuum $\boldsymbol
\varphi[vac]$. It is a tensorial field $\tilde{\boldsymbol\varphi}$ 
associated with the bundles $\uU\!M$ and $S\uuU\!M$ of the same type as
$\boldsymbol\varphi[vac]$:
$$
\hskip -2em
\tilde\varphi^{\kern 0.4pt\alpha111}=
\varphi^{\kern 0.4pt\alpha111}[vac]+\xi^{\kern 0.4pt\alpha111}.
\mytag{3.8}
$$
Note that $|\boldsymbol\varphi|\neq 0$. The perturbation $\boldsymbol
\xi$ in \mythetag{3.8} is assumed to be sufficiently small so that 
$|\tilde{\boldsymbol\varphi}|\neq 0$. By analogy to \mythetag{2.74}
we denote
$$
\hskip -2em
|\tilde{\boldsymbol\varphi}|^2=\frac{\tilde v^2}{2}\neq 0.
\mytag{3.9}
$$
Due to \mythetag{3.9} we can construct two tensor fields
$$
\xalignat 2
&\hskip -2em
\frac{\sqrt{2}\,\boldsymbol\varphi[vac]}{v},
&&\frac{\sqrt{2}\,\tilde{\boldsymbol\varphi}}{\tilde v}
\mytag{3.10}
\endxalignat
$$
and for each point $p\in M$ treat their values as two vectors of 
the unit length in a two-dimensional Hermitian vector space (the 
Hermitian form of this space is given by the formula \mythetag{2.64}).
\mylemma{3.1} For any two vectors $\bold v$ and $\tilde{\bold v}$ of 
the unit length in a Hermitian space of the dimension $2$ or higher 
there is a unitary operator $\boldsymbol\Omega$ with the unit determinant
$\det\boldsymbol\Omega=1$ such that $\tilde{\bold v}=\boldsymbol\Omega
(\bold v)$.
\endproclaim
    This lemma is a rather simple fact from the linear algebra. Its proof
is left to the reader. Applying this lemma to the vectors \mythetag{3.10},
we find that there is a special unitary matrix $\Omega\in\MatGrSU(2)$ and 
a gauge transformation \mythetag{3.3} such that
$$
\hskip -2em
\boldsymbol\varphi[vac]\longmapsto
\frac{v}{\tilde v}\,\tilde{\boldsymbol\varphi}.
\mytag{3.11}
$$
The formula \mythetag{3.11} means that each sufficiently small perturbation
of the Higgs vacuum is decomposed into elongation and rotation parts and the
rotation part is gauge equivalent to the initial vacuum. For this reason
only elongation type perturbations of the Higgs vacuum are considered in the
Standard Model (see \mycite{3}):
$$
\hskip -2em
\boldsymbol\varphi=\boldsymbol\varphi[vac]
+\frac{\chi}{v}\,\boldsymbol\varphi[vac].
\mytag{3.12}
$$
Here $\chi$ is a real scalar field. It is also called the {\it Higgs field}.
This terminology makes no confusion since both Higgs fields are closely 
related through the formula \mythetag{3.12}.
\head
4. Symmetry breaking\\
and boson fields recalculation.
\endhead
     From now on assume that some nontrivial Higgs vacuum is fixed. 
On fixing it we say that the initial $\MatGrSU(2)\times\MatGrU(1)$ 
symmetry is broken. However, some definite amount of this symmetry 
remains unbroken. Remember that $\nabla\boldsymbol\varphi=0$ with 
respect to the vacuum connections in \mythetag{2.36} and \mythetag{2.37}. 
Let's study those gauge fields $\bolduCA$ and $\bolduuCA$ in \mythetag{2.36} 
and \mythetag{2.37} that preserve this equality:
$$
\nabla_{\!k}\,\varphi^{\kern 0.4pt\alpha111}[vac]
-\frac{i\,e}{\hbar\,c}\left(\,\shave{\sum^2_{\theta=1}}
g_2\,\uuCA^\alpha_{k\kern 0.4pt\theta}
\,\varphi^{\kern 1.0pt\theta\kern 0.4pt 111}[vac]
+3\,g_1\,\uCA^1_{k1}\,\varphi^{\kern 0.6pt\alpha111}[vac]
\!\right)=0.\qquad
\mytag{4.1}
$$
Since $\nabla_{\!k}\,\varphi^{\kern 0.4pt\alpha111}[vac]=0$, the equality
\mythetag{4.1} reduces to a non-differential equality:
$$
\hskip -2em
\sum^2_{\theta=1}g_2\,\uuCA^\alpha_{k\kern 0.4pt\theta}
\,\varphi^{\kern 1.0pt\theta\kern 0.4pt 111}[vac]
+3\,g_1\,\uCA^1_{k1}\,\varphi^{\kern 0.6pt\alpha111}[vac]=0.
\mytag{4.2}
$$
This equality \mythetag{4.2} here should be treated as a complement 
to the concordance conditions \mythetag{2.48}, \mythetag{2.49}, and
\mythetag{2.50}. In order to solve the total set of the equations
\mythetag{4.2}, \mythetag{2.48}, \mythetag{2.49}, \mythetag{2.50} 
we introduce a projector $\bold P[vac]$ with the components:
$$
\hskip -2em
P^\beta_\alpha[vac]=\sum^2_{\bar a=1}\frac{\uuD_{\alpha\bar\alpha}
\,\uD_{11}\,\uD_{11}\,\uD_{11}}{|\boldsymbol\varphi[vac]|^2}
\,\overline{\varphi^{\kern 0.4pt\raise 1.2pt\hbox{$\ssize\bar\alpha
111$}}[vac]}\,\varphi^{\kern 0.4pt\beta 111}[vac].
\mytag{4.3}
$$
The projection operator with the components \mythetag{4.3} is an
orthogonal projector since it is a Hermitian operator in the sense of 
the following equality:
$$
\hskip -2em
\sum^2_{a=1}\uuD_{a\bar j}
\,P^{\kern 0.2pt a}_i[vac]=\sum^2_{\bar a=1}\uuD_{i\bar a}
\,\overline{P^{\kern 0.2pt\raise 0.6pt
\hbox{$\ssize\bar a$}}_{\bar j}[vac]}.
\mytag{4.4}
$$
Now, in addition to $\boldsymbol\varphi[vac]$, we introduce another
field $\boldsymbol\phi[vac]$ derived from it. This field can also be 
called the Higgs vacuum field:
$$
\hskip -2em
\phi^{\kern 0.8pt\alpha}_{111}[vac]=\sum^2_{\beta=1}\sum^2_{\bar\beta=1}
\uud^{\kern 0.4pt\alpha\beta}\,\uuD_{\beta\bar\beta}\,\uD_{11}\,\uD_{11}
\,\uD_{11}\,\overline{\varphi^{\kern 0.4pt\raise 1.2pt\hbox{$\ssize
\bar\beta 111$}}[vac]}.
\mytag{4.5}
$$
The modulus of the field \mythetag{4.5} is determined by the formula
similar to \mythetag{2.64}:
$$
\hskip -2em
|\boldsymbol\phi|^2=
\sum^2_{\alpha=1}\sum^2_{\bar\alpha=1}\uuD_{\alpha\bar\alpha}\,
\uD^{11}\,\uD^{11}\,\uD^{11}\,\phi^{\kern 0.8pt\alpha}_{111}\,
\overline{\phi^{\kern 0.8pt\raise 
1.2pt\hbox{$\ssize\bar\alpha$}}_{111}}.
\mytag{4.6}
$$
Using \mythetag{2.16}, \mythetag{2.32}, and \mythetag{2.33}, from 
\mythetag{2.64}, \mythetag{2.73}, and \mythetag{4.6} one easily 
derives 
$$
\hskip -2em
|\boldsymbol\phi[vac]|^2=|\boldsymbol\varphi[vac]|^2=\frac{v^2}{2}
=\const.
\mytag{4.7}
$$
Due to \mythetag{4.7}, using $\boldsymbol\phi[vac]$, one can define
another projection operator $\bold Q[vac]$:
$$
\hskip -2em
Q^\beta_\alpha[vac]=\sum^2_{\bar a=1}\frac{\uuD_{\alpha\bar\alpha}
\,\uD^{11}\,\uD^{11}\,\uD^{11}}{|\boldsymbol\phi[vac]|^2}
\,\overline{\phi^{\kern 0.4pt\raise 1.2pt\hbox{$\ssize\bar
\alpha$}}_{111}[vac]}\,\phi^{\kern 0.4pt\beta}_{111}[vac].
\mytag{4.8}
$$
Like the projector $\bold P[vac]$ with the components \mythetag{4.3},
the projector $\bold Q[vac]$ witn the components \mythetag{4.8} is an
orthogonal projector. This fact is expressed by the equality
$$
\hskip -2em
\sum^2_{a=1}\uuD_{a\bar j}
\,Q^{\kern 0.2pt a}_i[vac]=\sum^2_{\bar a=1}\uuD_{i\bar a}
\,\overline{Q^{\kern 0.2pt\raise 0.6pt
\hbox{$\ssize\bar a$}}_{\bar j}[vac]}.
\mytag{4.9}
$$
It easy to see that tne equality \mythetag{4.9} is analogous to
\mythetag{4.4}.\par
     Two Higgs fields $\boldsymbol\phi[vac]$ and $\boldsymbol\varphi
[vac]$ are perpendicular to each other in the sense of the following
equality reflecting the Hermitian forms in \mythetag{2.64} and
\mythetag{4.6}:
$$
\sum^2_{\alpha=1}\sum^2_{\bar\alpha=1}\uuD_{\alpha\bar\alpha}\,
\overline{\phi^{\kern 0.4pt\raise 1.2pt\hbox{$\ssize\bar
\alpha$}}_{111}[vac]}\,\varphi^{\kern 0.4pt\alpha111}[vac]=0.
\mytag{4.10}
$$
Due \mythetag{4.10} the projectors $\bold P[vac]$ and $\bold Q[vac]$ 
are complementary to each other:
$$
\hskip -2em
\bold P[vac]+\bold Q[vac]=\idop.
\mytag{4.11}
$$
In a coordinate form the operator equality \mythetag{4.11} looks like
$$
\hskip -2em
P^\beta_\alpha[vac]+Q^\beta_\alpha[vac]=\delta^\beta_\alpha.
\mytag{4.12}
$$
\par
    Apart from \mythetag{4.11} and \mythetag{4.12}, there are some other 
consequences of the orthogonality condition \mythetag{4.10}. Though
$\boldsymbol\phi[vac]$ and $\boldsymbol\varphi[vac]$ are two tensorial
fields of the different types, due to \mythetag{4.7} and \mythetag{4.10}
one can treat them as a kind of orthogonal frame in the two-dimensional
complex bundle $S\uuU\!M$. Relying upon this interpretation of $\boldsymbol
\phi[vac]$ and $\boldsymbol\varphi[vac]$, we introduce two tensor fields 
$\bold W[\varphi\!\blacktriangleright\!\phi]$ and $\bold W[\phi\!
\blacktriangleright\!\varphi]$:
$$
\align
&\hskip -2em
W^\beta_{\alpha111111}[\varphi\!\blacktriangleright\!\phi]
=\sum^2_{\bar a=1}\frac{\uuD_{\alpha\bar\alpha}\,\uD_{11}
\,\uD_{11}\,\uD_{11}}{|\boldsymbol\varphi[vac]|^2}
\,\overline{\varphi^{\kern 0.4pt\raise 1.2pt\hbox{$\ssize
\bar\alpha111$}}[vac]}\,\phi^{\kern 0.4pt\beta}_{111}[vac],
\qquad
\mytag{4.13}\\
&\hskip -2em
W^{\beta111111}_\alpha[\phi\!\blacktriangleright\!\varphi]
=\sum^2_{\bar a=1}\frac{\uuD_{\alpha\bar\alpha}\,\uD^{11}
\,\uD^{11}\,\uD^{11}}{|\boldsymbol\phi[vac]|^2}
\,\overline{\phi^{\kern 0.4pt\raise 1.2pt\hbox{$\ssize\bar
\alpha$}}_{111}[vac]}\,\varphi^{\kern 0.4pt\beta 111}[vac].
\qquad
\mytag{4.14}
\endalign
$$
It easy to derive the following identities for the tensor fields
\mythetag{4.13} and \mythetag{4.14}:
$$
\gather
\hskip -2em
\sum^2_{\alpha=1}W^\beta_{\alpha111111}[\varphi\!\blacktriangleright
\!\phi]\ W^\alpha_{\gamma111111}[\varphi\!\blacktriangleright\!\phi]=0,
\mytag{4.15}\\
\hskip -2em
\sum^2_{\alpha=1}W^{\beta111111}_\alpha[\phi\!\blacktriangleright
\!\varphi]\ W^{\alpha111111}_\gamma[\phi\!\blacktriangleright\!\varphi]=0,
\mytag{4.16}\\
\hskip -2em
\sum^2_{\alpha=1}W^\beta_{\alpha111111}[\varphi\!\blacktriangleright
\!\phi]\ W^{\alpha111111}_\gamma[\phi\!\blacktriangleright\!\varphi]
=Q^\beta_\gamma[vac],
\mytag{4.17}\\
\hskip -2em
\sum^2_{\alpha=1}W^{\beta111111}_\alpha[\phi\!\blacktriangleright
\!\varphi]\,W^\alpha_{\gamma111111}[\varphi\!\blacktriangleright
\!\phi]=P^\beta_\gamma[vac].
\mytag{4.18}
\endgather
$$
In a coordinate-free form the identties \mythetag{4.15}, \mythetag{4.16},
\mythetag{4.17}, \mythetag{4.18} are written as
$$
\xalignat 2
&\hskip -2em
\bold W[\varphi\!\blacktriangleright\!\phi]^2=0,
&&\bold W[\varphi\!\blacktriangleright\!\phi]
\compos\bold W[\phi\!\blacktriangleright\!\varphi]
=\bold Q[vac],\quad
\mytag{4.19}\\
&\hskip -2em
\bold W[\phi\!\blacktriangleright\!\varphi]^2=0,
&&\bold W[\phi\!\blacktriangleright\!\varphi]
\compos\bold W[\varphi\!\blacktriangleright\!\phi]
=\bold P[vac].\quad
\mytag{4.20}
\endxalignat
$$
In addition to \mythetag{4.19} and \mythetag{4.20}, we have
$$
\xalignat 2
&\hskip -2em
\bold W[\varphi\!\blacktriangleright\!\phi]\compos\bold P[vac]
=\bold W[\varphi\!\blacktriangleright\!\phi],
&&\bold W[\varphi\!\blacktriangleright\!\phi]
\compos\bold Q[vac]=0,
\quad
\mytag{4.21}\\
&\hskip -2em
\bold Q[vac]\compos\bold W[\varphi\!\blacktriangleright\!\phi]
=\bold W[\varphi\!\blacktriangleright\!\phi],
&&\bold P[vac]\compos\bold W[\varphi\!\blacktriangleright\!\phi]=0,
\quad
\mytag{4.22}\\
&\hskip -2em
\bold W[\phi\!\blacktriangleright\!\varphi]\compos\bold Q[vac]
=\bold W[\phi\!\blacktriangleright\!\varphi],
&&\bold W[\phi\!\blacktriangleright\!\varphi]\compos\bold P[vac]
=0,\quad
\mytag{4.23}\\
&\hskip -2em
\bold P[vac]\compos\bold W[\phi\!\blacktriangleright\!\varphi]
=\bold W[\phi\!\blacktriangleright\!\varphi],
&&\bold Q[vac]\compos\bold W[\phi\!\blacktriangleright\!\varphi]
=0.\quad
\mytag{4.24}
\endxalignat
$$
\par
     Let $A^{\ldots\,\beta\,\ldots}_{\ldots\,\alpha\,\ldots}$ be the
components of some arbitrary tensorial field $\bold A$ with at least 
one upper index and at least one lower index associated with the
two-dimensional bundle $S\uuU\!M$. By dots we denote other indices of
$\bold A$ that could be present. In this case we have the following
expansion for the tensor field $\bold A$:
$$
\hskip -2em
\gathered
A^{\ldots\,\beta\,\ldots}_{\ldots\,\alpha\,\ldots}
=A^{{\sssize +}\ldots}_{\ldots}\cdot P^\beta_\alpha[vac]
+A^{{\sssize -}\ldots}_{\ldots}\cdot Q^\beta_\alpha[vac]\,+\\
+\,W^{{\sssize +}111111\,\ldots}_{\ldots\,\ldots\,\ldots\,\ldots}
\cdot W^\beta_{\alpha111111}[\varphi\!\blacktriangleright\!\phi]
+W^{\sssize -\,\ldots\,\ldots\,\ldots}_{111111\,\ldots}\cdot 
W^{\beta111111}_\alpha[\phi\!\blacktriangleright\!\varphi].
\endgathered
\mytag{4.25}
$$
The expansion \mythetag{4.25} is derived with the use of the formula
\mythetag{4.11} and the above identities \mythetag{4.21}, \mythetag{4.22},
\mythetag{4.23}, \mythetag{4.24}. Let's apply this expansion to the
components of the tensor $\bolduuCA$ in \mythetag{2.37}. As a result we
get
$$
\hskip -2em
\gathered
\uuCA^\beta_{k\alpha}=A^{\sssize+}_k\cdot P^\beta_\alpha[vac]+
A^{\sssize-}_k\cdot Q^\beta_\alpha[vac]\,+\\
+\,W^{{\sssize +}111111}_k\cdot W^\beta_{\alpha111111}[\varphi
\!\blacktriangleright\!\phi]+W^{\sssize -}_{k111111}\cdot 
W^{\beta111111}_\alpha[\phi\!\blacktriangleright\!\varphi].
\endgathered
\mytag{4.26}
$$
Here $A^{\sssize+}_k$ and $A^{\sssize-}_k$ are the components of two
fields which are $\MatGrSU(2)$-singlets. Applying \mythetag{2.49} and
\mythetag{2.50} to \mythetag{4.26}, we find that 
$$
\xalignat 3
&\hskip -2em
A^{\sssize+}_k=\overline{A^{\sssize+}_k},
&&A^{\sssize-}_k=\overline{A^{\sssize-}_k},
&&A^{\sssize-}_k=-A^{\sssize+}_k.
\qquad
\mytag{4.27}
\endxalignat
$$
In other words $A^{\sssize+}_k$ and $A^{\sssize-}_k$ are the 
components of two mutually opposite covectorial fields. The index 
$k$ in \mythetag{4.26} and \mythetag{4.27} is a spatial index 
associated with the tangent bundle $TM$. In physical literature 
the following notation is a tradition (see \mycite{3}):
$$
\xalignat 3
&\hskip -2em
A^{\sssize+}_k=-A^3_k, 
&&A^{\sssize-}_k=A^3_k, 
&&A^3_k=\overline{A^{\raise 0.6pt\hbox{$\ssize 3$}}_k}.
\qquad
\mytag{4.28}
\endxalignat
$$
The number $3$ in \mythetag{4.28} is not a tensorial index. In \mycite{3}
and in many other books and papers it is set because $A^3_k$ first appear
as coefficients of the third Pauli matrix
$$
\boldsymbol\sigma_3=\Vmatrix 1 & 0\\0 & -1\endVmatrix.
$$
For $W^{{\sssize +}111111}_k$ and $W^{\sssize -}_{k111111}$ in 
\mythetag{4.26} from \mythetag{2.49} we derive
$$
\hskip -2em
\aligned
W^{\sssize -}_{k111111}=D_{11}\,D_{11}\,D_{11}\,D_{11}\,D_{11}\,D_{11}\,
\overline{W^{{\sssize +}111111}_k},\\
W^{{\sssize +}111111}_k=D^{11}\,D^{11}\,D^{11}\,D^{11}\,D^{11}\,D^{11}\,
\overline{W^{\sssize -}_{k111111}}.
\endaligned
\mytag{4.29}
$$
The index $k$ in \mythetag{4.29} is associated with the tangent 
bundle $TM$. For this reason $W^{{\sssize +}111111}_k$ and 
$W^{\sssize -}_{k111111}$ are the components of two covectorial 
fields $\bold W^{\sssize +}$ and $\bold W^{\sssize -}$. These
fields correspond to $W$-bosons. They are scalar fields with respect 
to the bundle $S\uuU\!M$, i\.\,e\. they are $\MatGrSU(2)$-singlets,
but they are not scalar with respect to the one-dimensional bundle 
$\uU\!M$.\par
    Thus the relationships \mythetag{2.49} and \mythetag{2.50} are
resolved to \mythetag{4.28} and  \mythetag{4.29}. Let's apply 
\mythetag{4.26} to \mythetag{4.2}. As a result, taking into account
\mythetag{4.28} and  \mythetag{4.29}, we get
$$
\xalignat 2
&\hskip -2em
W^{\sssize -}_{k111111}=0,
&&W^{{\sssize +}111111}_k=0,
\quad\\
\vspace{-1.0ex}
&&&\mytag{4.30}\\
\vspace{-2.0ex}
&\hskip -2em
A^{\sssize+}_k=\frac{3\,g_1}{g_1}\,\uCA^1_{k1},
&&A^{\sssize-}_k=-\frac{3\,g_1}{g_1}\,\uCA^1_{k1}.
\quad
\endxalignat
$$
The formulas \mythetag{4.30} show that there is a special combination 
of the gauge fields $\bolduCA$ and $\bolduuCA$ which, upon substituting 
into \mythetag{2.36} and \mythetag{2.37}, preserves the classical Higgs 
vacuum. The gauge field $\bolduuCA$ is expressed through the gauge
field $\bolduCA$ in this combination. In other words this fact means 
that the $4$-dimensional gauge symmetry $\MatGrU(2)\times\MatGrU(1)$
reduces to $1$-dimensional gauge symmetry $\MatGrU(1)$. However, this
new $\MatGrU(1)$ symmetry does not coincide with the initial $\MatGrU(1)$
symmetry being a component in the direct product $\MatGrSU(2)\times
\MatGrU(1)$. This residual $\MatGrU(1)$ symmetry in the electro-weak
bundles $S\uuU\!M$ and $\uU\!M$ is interpreted as a $\MatGrU(1)$
symmetry of the electromagnetism. In order to reveal this hidden 
$\MatGrU(1)_{\text{em}}$ symmetry the following two real fields $\bold A$ 
and $\bold Z$ are introduced in addition to $\bold W^{\sssize +}$ and
$\bold W^{\sssize -}$:
$$
\xalignat 2
&\hskip -2em
Z_k=\frac{-g_2\,A^3_k+3\,g_1\,\uCA^1_{k1}}
{\sqrt{(g_2)^2+(3\,g_1)^2}},
&&A_k=\frac{3\,g_1\,A^3_k+g_2\,\uCA^1_{k1}}
{\sqrt{(g_2)^2+(3\,g_1)^2}}.
\qquad
\mytag{4.31}
\endxalignat
$$
The relationships \mythetag{4.30} now mean that \mythetag{4.2} is fulfilled
provided 
$$
\xalignat 3
&\hskip -2em
\bold W^{\sssize +}=0,
&&\bold W^{\sssize -}=0,
&&\bold Z=0.
\quad
\mytag{4.32}
\endxalignat
$$
The condition $\bold A=0$ is absent in \mythetag{4.32} because the 
classical Higgs vacuum $\boldsymbol\varphi[vac]$ is invariant with 
respect to purely electromagnetic gauge transformations.\par
    Using the relationships \mythetag{4.31}, we can express $A^3_k$
and $\uCA^1_{k1}$ through $A_k$ and $Z_k$. For the component of the
gauge field $\bolduCA$ we derive:
$$
\hskip -2em
\uCA^1_{k1}=\frac{g_2\,A_k+3\,g_1\,Z_k}{\sqrt{(g_2)^2+(3\,g_1)^2}}.
\mytag{4.33}
$$
In a similar way, from \mythetag{4.31} for $A^3_k$ we derive
$$
\hskip -2em
A^3_k=\frac{3\,g_1\,A_k-g_2\,Z_k}{\sqrt{(g_2)^2+(3\,g_1)^2}}.
\mytag{4.34}
$$
Then we substitute \mythetag{4.34} back into \mythetag{4.28} and
\mythetag{4.26}. As a result we get
$$
\hskip -2em
\gathered
\uuCA^\beta_{k\alpha}=\frac{3\,g_1\,A_k-g_2\,Z_k}
{\sqrt{(g_2)^2+(3\,g_1)^2}}\cdot\bigl(Q^\beta_\alpha[vac]-
P^\beta_\alpha[vac]\,\bigr)\,+\\
\vspace{2ex}
+\,W^{{\sssize +}111111}_k\cdot W^\beta_{\alpha111111}[\varphi
\!\blacktriangleright\!\phi]+W^{\sssize -}_{k111111}\cdot 
W^{\beta111111}_\alpha[\phi\!\blacktriangleright\!\varphi].
\endgathered
\mytag{4.35}
$$
The next step is to substitute \mythetag{4.35} into \mythetag{2.45}.
We do it in several substeps:
$$
\gather
\hskip -2em
\sum^2_{h=1}\uuCA^p_{i\,h}\,\uuCA^{\!h}_{j\,k}
=A^3_i\,A^3_j\cdot\delta^p_q+\bigl(A^3_i\,W^{{\sssize +}111111}_j
-A^3_j\,W^{{\sssize +}111111}_i\bigr)\,\cdot\\
\hskip -2em
\cdot\,W^p_{q111111}[\varphi\!\blacktriangleright\!\phi]
+\bigl(A^3_j\,W^{\sssize -}_{i111111}-A^3_i
\,W^{\sssize -}_{j111111}\bigr)\cdot W^{p111111}_q[\phi
\!\blacktriangleright\!\varphi]\,+
\mytag{4.36}\\
\vspace{2ex}
\hskip -2em
+\,W^{{\sssize +}111111}_i\,W^{\sssize -}_{j111111}
\cdot Q^p_q[vac]+W^{{\sssize +}111111}_j\,W^{\sssize -}_{i111111}
\cdot P^p_q[vac].
\endgather
$$
From \mythetag{4.36} by means of alternation we immediately derive
$$
\gather
\hskip -2em
\sum^2_{h=1}\left(\uuCA^p_{i\,h}\,\uuCA^{\!h}_{j\,q}-\uuCA^p_{j\,h}
\,\uuCA^{\!h}_{i\,q}\right)
=2\,\bigl(A^3_i\,W^{{\sssize +}111111}_j
-A^3_j\,W^{{\sssize +}111111}_i\bigr)\,\cdot
\quad\\
\hskip -2em
\cdot\,W^p_{q111111}[\varphi\!\blacktriangleright\!\phi]
+2\,\bigl(A^3_j\,W^{\sssize -}_{i111111}-A^3_i
\,W^{\sssize -}_{j111111}\bigr)\cdot W^{p111111}_q[\phi
\!\blacktriangleright\!\varphi]\,+\quad
\mytag{4.37}\\
\vspace{2ex}
\hskip -2em
+\,\left(W^{{\sssize +}111111}_i\,W^{\sssize -}_{j111111}
-W^{{\sssize +}111111}_j\,W^{\sssize -}_{i111111}\right)\cdot 
\left(Q^p_q[vac]-P^p_q[vac]\right).\quad
\endgather
$$
For the differential part of the tensor $\bolduuF$ in the formula 
\mythetag{2.45} we get
$$
\gather
\hskip -2em
\nabla_{\!i}\,\uuCA^p_{j\,q}-\nabla_{\!j}\,\uuCA^p_{i\,q}
=\frac{Q^p_q[vac]-P^p_q[vac]}{\sqrt{(g_2)^2+(3\,g_1)^2}}
\cdot\bigl(3\,g_1\,(\nabla_{\!i}A_j-\nabla_{\!j}A_i)\,-\quad\\
\vspace{2ex}
\hskip -2em
-\,g_2\,(\nabla_{\!i}Z_j-\nabla_{\!j}Z_i)\bigr)+
(\nabla_{\!i}W^{{\sssize +}111111}_j-\nabla_{\!j}
W^{{\sssize +}111111}_i)\,\cdot\quad
\mytag{4.38}\\
\vspace{2ex}
\hskip -2em
\cdot\,W^p_{q111111}[\varphi\!\blacktriangleright\!\phi]
+(\nabla_{\!i}W^{\sssize -}_{j111111}-\nabla_{\!j}
W^{\sssize -}_{i\kern 0.2pt 111111})\cdot W^{p111111}_q[\phi
\!\blacktriangleright\!\varphi].\quad
\endgather
$$
Now let's introduce the following notations analogous to 
\mythetag{2.44}:
$$
\xalignat 2
&\hskip -2em
F_{ij}=\nabla_{\!i}A_j-\nabla_{\!j}A_i,
&&\Cal Z_{ij}=\nabla_{\!i}Z_j-\nabla_{\!j}Z_i.
\mytag{4.39}
\endxalignat
$$
In the case of the fields $\bold W^{\sssize +}$ and $\bold W^{\sssize -}$
we need more complicated notations:
$$
\align
&\hskip -2em
\aligned
\Cal W^{{\sssize +}111111}_{ij}
&=\nabla_{\!i}W^{{\sssize +}111111}_j-\nabla_{\!j}
W^{{\sssize +}111111}_i\,-\\
&-\,\frac{6\,i\,e\,g_1}{\hbar\,c}
\left(\uCA^1_{i\kern 0.2pt 1}\,W^{{\sssize +}111111}_j-\uCA^1_{j1}
\,W^{{\sssize +}111111}_i\right).
\endaligned
\mytag{4.40}\\
\vspace{2ex}
&\hskip -2em
\aligned
\Cal W^{\sssize -}_{ij111111}
&=\nabla_{\!i}W^{\sssize -}_{j111111}-\nabla_{\!j}
W^{\sssize -}_{i\kern 0.2pt 111111}\,+\\
&+\,\frac{6\,i\,e\,g_1}{\hbar\,c}
\left(\uCA^1_{i\kern 0.2pt 1}\,W^{\sssize -}_{j111111}
-\uCA^1_{j1}\,W^{\sssize -}_{i\kern 0.2pt 111111}\right).
\endaligned
\mytag{4.41}
\endalign
$$
Here $\uCA^1_{i\kern 0.2pt 1}$ and $\uCA^1_{j1}$ are determined by 
the formula \mythetag{4.33}. Now we need to apply \mythetag{4.39}, 
\mythetag{4.40}, and \mythetag{4.41} to \mythetag{4.38}, and combine
\mythetag{4.38} with \mythetag{4.37} according to the formula
\mythetag{2.45}. As a result we obtain
$$
\hskip -2em
\gathered
\uuF^p_{\!\!qij}=\left(Q^p_q[vac]-P^p_q[vac]\right)\cdot\!\left(
\frac{3\,g_1\,F_{ij}-\,g_2\,\Cal Z_{ij}}{\sqrt{(g_2)^2+(3\,g_1)^2}}
\,-\right.\\
\left.-\frac{i\,e\,g_2}{\hbar\,c}\left(W^{{\sssize +}111111}_i
\,W^{\sssize -}_{j111111}-W^{{\sssize +}111111}_j
\,W^{\sssize -}_{i111111}\right)\!\!\vphantom{\frac{3\,g_1\,F_{ij}-\,g_2
\,\Cal Z_{ij}}{\sqrt{(g_2)^2+(3\,g_1)^2}}}\right)+\\
+\,W^p_{q111111}[\varphi\!\blacktriangleright\!\phi]
\cdot\!\left(\Cal W^{{\sssize +}111111}_{ij}
+\frac{2\,i\,e\,\sqrt{(g_2)^2+(3\,g_1)^2}}{\hbar\,c}\,\times\right.\\
\left.\times\,\left(Z_i\,W^{{\sssize +}111111}_j-Z_j
\,W^{{\sssize +}111111}_i\right)\!\!\vphantom{\frac{2\,i\,e
\,\sqrt{(g_2)^2+(3\,g_1)^2}}{\hbar\,c}}\right)+\\
+\,W^{p111111}_q[\phi\!\blacktriangleright\!\varphi]
\cdot\!\left(\Cal W^{\sssize -}_{ij111111}
-\frac{2\,i\,e\,\sqrt{(g_2)^2+(3\,g_1)^2}}{\hbar\,c}\,\times\right.\\
\left.\times\,\left(Z_i\,W^{\sssize -}_{j111111}
-Z_j\,W^{\sssize -}_{i\kern 0.2pt 111111}\right)\!\!\vphantom{\frac{2
\,i\,e\,\sqrt{(g_2)^2+(3\,g_1)^2}}{\hbar\,c}}
\right).
\endgathered
\mytag{4.42}
$$\par
    Thus, in \mythetag{4.42} we have expressed the field $\bolduuF$
through the fields $\bold W^{\sssize +}$, $\bold W^{\sssize -}$, and
$\bold Z$ corresponding to $W$ and $Z$-bosons of the Standard Model
and through the field $\bold A$ which is interpreted as the 
$4$-dimensional vector-potential of the electromagnetic field. Our
next goal is to express $\bolduF$ through $\bold Z$ and $\bold A$
using the formulas \mythetag{4.33} and  \mythetag{2.44}. Substituting
\mythetag{4.33} into \mythetag{2.44} and taking into account
\mythetag{4.39}, we derive
$$
\hskip -2em
\uF_{\!\!ij}=\frac{g_2\,F_{ij}+3\,g_1\,\Cal Z_{ij}}
{\sqrt{(g_2)^2+(3\,g_1)^2}}\,.
\mytag{4.43}
$$
Then we substitute \mythetag{4.43} and \mythetag{4.42} into \mythetag{2.59} 
and \mythetag{2.60} respectively. As a result for the sum of two integrals
\mythetag{2.59} and \mythetag{2.60} we get
$$
\hskip -2em
\Cal L_1+\Cal L_2=\Cal L_{11}+\Cal L_{12}+\Cal L_{21}+\Cal L_{22}
+\Cal L_{23}+\Cal L_{24}+\Cal L_{25}+\Cal L_{26},
\mytag{4.44}
$$
where $\Cal L_{11}$ and $\Cal L_{12}$ are standard kinetic terms for
two real vector fields $\bold A$ and $\bold Z$:
$$
\align
\hskip -2em
\Cal L_{11}&=-\frac{1}{16\,\pi\,c}
\int\sum^3_{i=0}\sum^3_{j=0}F_{ij}\,F^{ij}\,dV,
\mytag{4.45}\\
\hskip -2em
\Cal L_{12}&=-\frac{1}{16\,\pi\,c}
\int\sum^3_{i=0}\sum^3_{j=0}\Cal Z_{ij}\,\Cal Z^{ij}\,dV.
\mytag{4.46}
\endalign
$$
The term $\Cal L_{21}$ in \mythetag{4.44} is a standard kinetic term 
for two mutually conjugate complex vector fields $\bold W^{\sssize +}$ 
and $\bold W^{\sssize -}$ (see \mythetag{4.29}):
$$
\hskip -2em
\Cal L_{21}=-\frac{1}{16\,\pi\,c}
\int\sum^3_{i=0}\sum^3_{j=0}\Cal W^{{\sssize +}111111}_{ij}
\,\Cal W^{{\sssize -}ij}_{111111}\,dV.
\mytag{4.47}
$$
The term $\Cal L_{22}$ is a purely potential term responsible for the
self-action of $W$-bosons:
$$
\gathered
\Cal L_{22}=-\frac{\left(\dfrac{e\,g_2}
{\hbar\,c}\right)^{\!2}}{8\,\pi\,c}\int\sum^3_{i=0}\sum^3_{j=0}
W^{{\sssize +}111111}_i\,W^{{\sssize -}i}_{111111}\,
W^{{\sssize +}111111}_j\,W^{{\sssize -}j}_{111111}\,dV\,+\\
+\,\frac{\left(\dfrac{e\,g_2}
{\hbar\,c}\right)^{\!2}}{8\,\pi\,c}\int\sum^3_{i=0}\sum^3_{j=0}
W^{{\sssize +}111111}_i\,W^{{\sssize +}i\kern 0.6pt 111111}
\,W^{\sssize -}_{j\kern 0.6pt 111111}\,W^{{\sssize -}j}_{111111}\,dV.
\endgathered
\mytag{4.48}
$$
The term $\Cal L_{23}$ is a mixed term responsible the interaction
of $W$-bosons and photons:
$$
\hskip -2em
\Cal L_{23}=\frac{\dfrac{3\,i\,e\,g_2\,g_1}{8\,\pi\,\hbar\,c^2}}
{\sqrt{(g_2)^2+(3\,g_1)^2}}\int\sum^3_{i=0}\sum^3_{j=0}F^{ij}
\,W^{{\sssize +}111111}_i\,W^{\sssize -}_{j\kern 0.6pt 111111}\,dV.
\mytag{4.49}
$$
The term $\Cal L_{24}$ is a mixed term responsible the interaction
of $W$ and $Z$-bosons:
$$
\hskip -2em
\Cal L_{24}=-\frac{\dfrac{i\,e\,(g_2)^2}{8\,\pi\,\hbar\,c^2}}
{\sqrt{(g_2)^2+(3\,g_1)^2}}\int\sum^3_{i=0}\sum^3_{j=0}\Cal Z^{ij}
\,W^{{\sssize +}111111}_i\,W^{\sssize -}_{j\kern 0.6pt 111111}\,dV.
\mytag{4.50}
$$
The term $\Cal L_{25}$ is a potential term responsible for the
interaction of $W$ and $Z$-bosons:
$$
\gathered
\Cal L_{25}=\frac{e^2((g_2)^2+(3\,g_1)^2)}
{2\,\pi\,\hbar^2\,c^3}\int\sum^3_{i=0}\sum^3_{j=0}
Z_i\,Z^i\,W^{{\sssize +}111111}_j
\,W^{{\sssize -}j}_{111111}\,dV\,-\\
-\,\frac{e^2((g_2)^2+(3\,g_1)^2)}
{2\,\pi\,\hbar^2\,c^3}\int\sum^3_{i=0}\sum^3_{j=0}
Z^i\,W^{{\sssize +}111111}_i\,Z_j\,W^{{\sssize -}j}_{111111}
\,dV.
\endgathered
\mytag{4.51}
$$
The term $\Cal L_{26}$ is a mixed term responsible the interaction
of $W$ and $Z$-bosons:
$$
\gathered
\Cal L_{25}=\frac{i\,e\sqrt{(g_2)^2+(3\,g_1)^2}}
{4\,\pi\,\hbar\,c^2}\int\sum^3_{i=0}\sum^3_{j=0}
\Cal W^{{\sssize +}111111}_{ij}\,Z^i
\,W^{{\sssize -}j}_{111111}\,dV\,-\\
-\,\frac{i\,e\sqrt{(g_2)^2+(3\,g_1)^2}}
{4\,\pi\,\hbar\,c^2}\int\sum^3_{i=0}\sum^3_{j=0}
\Cal W^{{\sssize -}ij}_{111111}\,Z_i
\,W^{{\sssize +}111111}_j\,dV.
\endgathered
\mytag{4.52}
$$
The terms \mythetag{4.49}, \mythetag{4.50}, \mythetag{4.51}, and
\mythetag{4.52} are cubic with respect to fields. The terms 
\mythetag{4.48} is the fourth order term. None of these terms 
can be treated as a mass term corresponding to the kinetic terms
\mythetag{4.45}, \mythetag{4.46}, or \mythetag{4.46}.\par
\head
5. Masses of boson fields.
\endhead
     In order to find the massive terms corresponding to the 
kinetic terms \mythetag{4.45}, \mythetag{4.46}, and \mythetag{4.47} 
we consider again the elongation type perturbation of the Higgs
vacuum \mythetag{3.12}. In a coordinate form it is given by the
formula 
$$
\hskip -2em
\varphi^{\kern 0.4pt\alpha111}
=\varphi^{\kern 0.4pt\alpha111}[vac]
+\frac{\chi}{v}\,\varphi^{\kern 0.4pt\alpha111}[vac].
\mytag{5.1}
$$
Having non-vacuum Higgs field \mythetag{5.1}, we should apply
non-vacuum covariant differentiation to it in \mythetag{2.67},
i\.\,e\. we should use non-vacuum connections \mythetag{2.36}
and \mythetag{2.37} instead of purely vacuum ones. Then we get
$$
\hskip -2em
\gathered
\nabla_{\!k}\varphi^{\kern 0.4pt\alpha111}
=-\frac{i\,e}{\hbar\,c}\left(\,\shave{\sum^2_{\theta=1}}
g_2\,\uuCA^\alpha_{k\kern 0.4pt\theta}\,\varphi^{\kern 1.0pt
\theta\kern 0.4pt 111}[vac]+3\,g_1\,\uCA^1_{k1}
\,\varphi^{\kern 0.6pt\alpha111}[vac]\!\right)\times\\
\vspace{2ex}
\times\left(1+\frac{\chi}{v}\right)+\frac{\nabla_k\chi}{v}
\,\varphi^{\kern 0.6pt\alpha111}[vac].
\endgathered
\mytag{5.2}
$$
Note that $\nabla_{\!k}$ in the left hand side of \mythetag{5.2} 
is a non-vacuum nabla, while $\nabla_k\chi$ is a vacuum covariant
derivative applied to the real scalar Higgs field $\chi$. Let's 
apply \mythetag{4.26} to $\uuCA^\alpha_{k\kern 0.4pt\theta}$ in
\mythetag{5.2}. As a result we derive
$$
\hskip -2em
\gathered
\nabla_{\!k}\varphi^{\kern 0.4pt\alpha111}
=-\frac{i\,e}{\hbar\,c}\left(g_2\,A^{\sssize+}_k+
3\,g_1\,\uCA^1_{k1}\right)\left(1+\frac{\chi}{v}\right)
\varphi^{\kern 0.6pt\alpha111}[vac]-\\
\vspace{1ex}
-\frac{i\,e}{\hbar\,c}\ g_2\,W^{{\sssize +}111111}_k
\left(1+\frac{\chi}{v}\right)\phi^{\kern 0.4pt\alpha}_{111}[vac]
+\frac{\nabla_k\chi}{v}\,\varphi^{\kern 0.6pt\alpha111}[vac].
\endgathered
\mytag{5.3}
$$
Now remember the formulas \mythetag{4.28} and \mythetag{4.31}. Then
\mythetag{5.3} is rewritten as 
$$
\hskip -2em
\gathered
\nabla_{\!k}\varphi^{\kern 0.4pt\alpha111}
=-\frac{i\,e}{\hbar\,c}\,\sqrt{(g_2)^2+(3\,g_1)^2}\ Z_k
\left(1+\frac{\chi}{v}\right)
\varphi^{\kern 0.6pt\alpha111}[vac]-\\
\vspace{1ex}
-\frac{i\,e}{\hbar\,c}\ g_2\,W^{{\sssize +}111111}_k
\left(1+\frac{\chi}{v}\right)\phi^{\kern 0.4pt\alpha}_{111}[vac]
+\frac{\nabla_k\chi}{v}\,\varphi^{\kern 0.6pt\alpha111}[vac].
\endgathered
\mytag{5.4}
$$
Let's substitute \mythetag{5.4} into the formula \mythetag{2.66} 
and take into account \mythetag{4.7} and the orthogonality condition
\mythetag{4.10}. This yields
$$
\hskip -2em
\gathered
|\nabla\boldsymbol\varphi|^2
=\frac{e^2\,((g_2)^2+(3\,g_1)^2)}{\hbar^2\,c^2}\,
\frac{(v+\chi)^2}{2}\sum^3_{i=0}Z_i\,Z^i\,+\\
+\,\frac{1}{2}\sum^3_{i=0}\sum^3_{j=0}g^{ij}\,\nabla_{\!i}\,\chi
\ \nabla_{\!j}\,\chi\,+\frac{e^2\,(g_2)^2}{\hbar^2\,c^2}
\,\frac{(v+\chi)^2}{2}\,\sum^3_{i=0}W^{{\sssize +}111111}_i
\,W^{{\sssize -}i}_{111111}.
\endgathered
\mytag{5.5}
$$
The next step is to substitute \mythetag{5.1} into \mythetag{2.64}.
As a result we get
$$
\hskip -2em
|\boldsymbol\varphi|^2=\frac{(v+\chi)^2}{2}.
\mytag{5.6}
$$
Then we substitute \mythetag{5.6} into \mythetag{2.65} and derive
$$
V(\boldsymbol\varphi)=\lambda\,\frac{(v+\chi)^4}{4}
-\mu^2\,\frac{(v+\chi)^2}{2}.
\mytag{5.7}
$$
It is clear that $V(\boldsymbol\varphi)$ in \mythetag{5.7} is a fourth
order polynomial with respect to $\chi$:
$$
\hskip -2em
\gathered
V(\boldsymbol\varphi)=\frac{\lambda}{4}\,\chi^4+\lambda\,v\,\chi^3
+\left(\frac{3\,\lambda}{2}\,v^2-\frac{\mu^2}{2}\right)\,\chi^2\,+\\
+\left(\lambda\,v^3-\mu^2\,v\right)\,\chi+\left(\frac{\lambda}{4}\,
v^4-\frac{\mu^2}{2}\,v^2\right).
\endgathered
\mytag{5.8}
$$
Due to the second equality \mythetag{2.74} the term linear in $\chi$
in the above polynomial \mythetag{5.8} does vanish. The polynomial
$V(\boldsymbol\varphi)$ simplifies to
$$
\hskip -2em
V(\boldsymbol\varphi)=\frac{\lambda}{4}\,\chi^4+\lambda\,v\,\chi^3
+\mu^2\,\chi^2-\frac{\mu^2}{4}\,v^2.
\mytag{5.9}
$$
Before substituting \mythetag{5.5} and \mythetag{5.9} into \mythetag{2.67}
and \mythetag{2.68}, let's denote
$$
\hskip -2em
m_\chi=2\,m_\varphi.
\mytag{5.10}
$$
Then, substituting \mythetag{5.5} and \mythetag{5.9} into \mythetag{2.67}
and \mythetag{2.68} and using \mythetag{5.10}, we get
$$
\hskip -2em
\Cal L_4+\Cal L_5=\Cal L_{41}+\Cal L_{42}+\Cal L_{43}+\Cal L_{44}
+\Cal L_{45}+\Cal L_{51}+\Cal L_{52}+\Cal L_{53}.
\mytag{5.11}
$$
The terms $\Cal L_{41}$ and $\Cal L_{51}$ in the sum \mythetag{5.11} 
\pagebreak are the kinetic term and the mass term of the real scalar 
Higgs field $\chi$ respectively. They are given by the formulas
$$
\allowdisplaybreaks
\align
&\hskip -2em
\Cal L_{41}=\frac{\hbar^2}{2\,m_\chi\,c}\int\sum^3_{i=0}\sum^3_{j=0}g^{ij}
\,\nabla_{\!i}\,\chi\ \nabla_{\!j}\,\chi\,dV,
\mytag{5.12}\\
\vspace{2ex}
&\hskip -2em
\Cal L_{51}=-\frac{m_\chi\,c}{2}\int\chi^2\,dV.
\mytag{5.13}
\endalign
$$
In order to fit \mythetag{5.12} and \mythetag{5.13} the parameter $\mu$ 
in \mythetag{5.9} should be chosen so that
$$
\mu^2=2.
$$
The term $\Cal L_{42}$ in \mythetag{5.11} is the mass term corresponding 
to the kinetic term \mythetag{4.46}:
$$
\hskip -2em
\Cal L_{42}=\frac{c\,m_Z^{\,2}}{8\,\pi\,\hbar^2}\int\sum^3_{i=0}
Z_i\,Z^i\,dV.
\mytag{5.14}
$$
Comparing \mythetag{5.14} and \mythetag{5.5}, we derive the formula for
the mass of $Z$-bosons:
$$
\hskip -2em
m_Z=\sqrt{\frac{4\,\pi\,((g_2)^2+(3\,g_1)^2)}{m_\chi}}
\,\frac{e\,v\,\hbar}{c^2}.
\mytag{5.15}
$$
The term $\Cal L_{43}$ in \mythetag{5.11} is the mass term corresponding 
to the kinetic term \mythetag{4.47}:
$$
\hskip -2em
\Cal L_{43}=\frac{c\,m_W^{\,2}}{8\,\pi\,\hbar^2}\int\sum^3_{i=0}
W^{{\sssize +}111111}_i\,W^{{\sssize -}i}_{111111}\,dV.
\mytag{5.16}
$$
Comparing \mythetag{5.16} and \mythetag{5.5}, we derive the formula for
the mass of $W$-bosons:
$$
\hskip -2em
m_W=\sqrt{\frac{4\,\pi\,(g_2)^2}{m_\chi}}
\,\frac{e\,v\,\hbar}{c^2}.
\mytag{5.17}
$$
The term $\Cal L_{44}$ in \mythetag{5.11} describes the interaction of
$Z$-bosons with the real scalar Higgs field. It is represented by the
following two integrals:
$$
\Cal L_{44}=\frac{c\,m_Z^{\,2}}{4\,\pi\,\hbar^2\,v}\int\sum^3_{i=0}
\chi\,Z_i\,Z^i\,dV
+\frac{c\,m_Z^{\,2}}{8\,\pi\,\hbar^2\,v^2}\int\sum^3_{i=0}
\chi^2\,Z_i\,Z^i\,dV.\qquad
\mytag{5.18}
$$
In a similar way, the term $\Cal L_{45}$ in the sum \mythetag{5.11}
describes the interaction of $W$-bosons with the real scalar Higgs 
field $\chi$:
$$
\pagebreak
\hskip -2em
\gathered
\Cal L_{45}=\frac{c\,m_W^{\,2}}{4\,\pi\,\hbar^2\,v}\int\sum^3_{i=0}
\chi\,W^{{\sssize +}111111}_i\,W^{{\sssize -}i}_{111111}\,dV\,+\\
+\,\frac{c\,m_W^{\,2}}{8\,\pi\,\hbar^2\,v^2}\int\sum^3_{i=0}
\chi^2\,W^{{\sssize +}111111}_i\,W^{{\sssize -}i}_{111111}\,dV.
\endgathered
\mytag{5.19}
$$
The term $\Cal L_{52}$ in the sum \mythetag{5.11} describes the
self-action of the Higgs field:
$$
\hskip -2em
\Cal L_{52}=-\frac{m_\chi\,c}{2\,v}\int\chi^3\,dV
-\frac{m_\chi\,c}{2\,v^2}\int\chi^4\,dV.
\mytag{5.20}
$$
The last term $\Cal L_{53}$ in the sum \mythetag{5.11} is a constant term:
$$
\hskip -2em
\Cal L_{53}=\frac{m_\chi\,c\,v^2}{4}\int dV.
\mytag{5.21}
$$
Usually constant terms like \mythetag{5.21} are omitted. However, we prefer
to keep it in \mythetag{5.11}, since in a non-flat space-time manifold $M$
it can contribute to the cosmological constant $\Lambda$ (see \mycite{8}).
\par
     Unlike the sum \mythetag{4.44}, the interaction and self-action 
terms \mythetag{5.18}, \mythetag{5.19}, and \mythetag{5.20} in the sum 
\mythetag{5.11} are purely potential terms. They do not comprise the
covariant derivatives of the fields $\bold Z$, $\bold W^{\sssize +}$,
$\bold W^{\sssize -}$, and $\chi$.
\head
6. The lepton masses.
\endhead
     Having derived the formulas \mythetag{5.15} and \mythetag{5.17} 
for the masses of $Z$ and $W$-bosons, now we return back to the leptons
represented in the table \mythetag{1.1}. The kinetic term for leptons
in the total action of the Standard Model is written as follows:
$$
\hskip -2em
\aligned
\Cal L_6&=i\,\hbar\int\!\!\sum_{i=e,\mu,\tau}
\sum^4_{a=1}\sum^4_{\bar a=1}\sum^4_{b=1}\sum^3_{q=0}
\sum^2_{\alpha=1}\sum^2_{\bar\alpha=1}
\uD^{11}\,\uD^{11}\,\uD^{11}\,\times\\
\vspace{2ex}
&\qquad\qquad\times\,\uuD_{\alpha\bar\alpha}\,
\overline{\bulpsi^{\raise 0.6pt \hbox{$\ssize\bar a\bar
\alpha$}}_{111}[i]}\,D_{a\bar a}\,\gamma^{\kern 0.5pt aq}_b
\ \nabla_{\!q}\bulpsi^{b\alpha}_{111}[i]\,dV\,+\\
\vspace{2ex}
&+\,i\,\hbar\int\!\!\sum_{i=e,\mu,\tau}
\sum^4_{a=1}\sum^4_{\bar a=1}\sum^4_{b=1}\sum^3_{q=0}
\sum^2_{\alpha=1}\sum^2_{\bar\alpha=1}
\uD^{11}\,\uD^{11}\,\uD^{11}\,\times\\
\vspace{2ex}
&\qquad\qquad\times\,\uD^{11}\,\uD^{11}\,\uD^{11}\,
\overline{\circpsi^{\raise 0.6pt \hbox{$\ssize
\bar a$}}_{111111}[i]}\,D_{a\bar a}\,\gamma^{\kern 0.5pt aq}_b
\ \nabla_{\!q}\circpsi^b_{111111}[i]\,dV.
\endaligned
\mytag{6.1}
$$ 
Mass terms corresponding to \mythetag{6.1} should yield masses for
the electron, muon, and tauon, but their neutrinos should remain
massless. For this reason mass terms for leptons are introduced 
through the interaction with the Higgs field:
$$
\hskip -2em
\aligned
\Cal L_7&=-\!\!\sum_{i=e,\mu,\tau}\!\!\!h[i]\!\int\!\sum^4_{a=1}
\sum^4_{\bar a=1}\sum^2_{\alpha=1}\sum^2_{\bar\alpha=1}
\uD^{11}\,\uD^{11}\,\uD^{11}\,\times\\
\vspace{2ex}
&\qquad\qquad\times\,\uuD_{\alpha\bar\alpha}\,D_{a\bar a}
\,\overline{\circpsi^{\raise 0.6pt \hbox{$\ssize\bar a$}}_{111111}[i]}
\ \overline{\varphi^{\kern 0.4pt\raise 1.2pt\hbox{$\ssize\bar\alpha111$}}}
\,\bulpsi^{a\alpha}_{111}[i]\,dV\,-\\
\vspace{2ex}
&-\!\sum_{i=e,\mu,\tau}\!\!\!h[i]\!\int\!\sum^4_{a=1}
\sum^4_{\bar a=1}\sum^2_{\alpha=1}\sum^2_{\bar\alpha=1}
\uD^{11}\,\uD^{11}\,\uD^{11}\,\times\\
\vspace{2ex}
&\qquad\qquad\times\,\uuD_{\alpha\bar\alpha}\,D_{a\bar a}
\,\overline{\bulpsi^{\raise 0.6pt \hbox{$\ssize\bar a
\bar\alpha$}}_{111}[i]}\ \varphi^{\alpha111}\,
\circpsi^a_{111111}[i]\,dV.
\endaligned
\mytag{6.2}
$$
Note that in \mythetag{6.1} and \mythetag{6.2} we used both doublet and
singlet wave functions \mythetag{1.4} and  \mythetag{1.5}. By $h[i]$ in
\mythetag{6.2} we denote three real constants $h[e]$, $h[\mu]$, $h[\tau]$
specific for each lepton generation.\par
     Now remember the projection operators $\bold P$ and $\bold Q$
introduced by the formulas \mythetag{4.3} and \mythetag{4.8}. Applying
them to doublet parts of lepton wave functions and taking into account
the equality \mythetag{4.12}, we get the expansion
$$
\hskip -2em
\bulpsi^{a\alpha}_{111}[i]=\bulpsi^a_{111111}[i]
\cdot\frac{\varphi^{\kern 0.4pt\alpha111}[vac]}{|\boldsymbol
\varphi[vac]|}+
\bulpsi^a[i]\cdot\frac{\phi^{\kern 0.8pt\alpha}_{111}[vac]}
{|\boldsymbol\phi[vac]|}.
\mytag{6.3}
$$
Note that $\bulpsi^a_{111111}[i]$ in \mythetag{6.3} are similar to 
$\circpsi^a_{111111}[i]$ in \mythetag{6.1} and \mythetag{6.2}. These
functions are chiral and antichiral parts of the complete wave functions:
$$
\hskip -2em
\psi^a_{111111}[i]=\bulpsi^a_{111111}[i]+\circpsi^a_{111111}[i].
\mytag{6.4}
$$
The singlet wave functions with the components \mythetag{6.4} describe
charged leptons: an electron for $i=e$, a muon for $i=\mu$, and a tauon
for $i=\tau$.\par
     The quantities $\bulpsi^a[i]$ are the components of other three 
singlet wave functions. They describe neutral leptons: a $e$-neutrino
for $i=e$, a $\mu$-neutrino for $i=\mu$, and a $\tau$-neutrino for 
$i=\tau$. These wave functions are chiral because antichiral (right)
neutrinos are not considered in the Standard Model.\par
{\bf A remark}. Charged leptons are distinguished from their neutral
counterparts due to the nontrivial Higgs vacuum that breaks the initial
$\MatGrSU(2)\times\MatGrU(1)$ symmetry.\par
     The covariant derivatives $\nabla_q$ in \mythetag{6.1} are complete
non-vacuum covariant derivatives. They are evaluated according to the
formulas \mythetag{2.3} and \mythetag{2.4}, where the electro-weak
connection components are taken from \mythetag{2.36} and \mythetag{2.37}
including the components of the gauge fields $\bolduCA$ and $\bolduuCA$.
Passing from non-vacuum to vacuum covariant derivatives in the formulas
\mythetag{2.3} and \mythetag{2.4}, we get
$$
\align
&\nabla_{\!q}\bulpsi^{b\kern 0.4pt\alpha}_{111}[i]
\to\nabla_{\!q}\bulpsi^{b\kern 0.4pt\alpha}_{111}[i]
-\frac{i\,e\,g_2}{\hbar\,c}\sum^2_{\theta=1}
\uuCA^\alpha_{q\kern 0.4pt\theta}
\,\bulpsi^{b\kern 0.8pt\theta}_{111}[i]
+\frac{3\,i\,e\,g_1}{\hbar\,c}\,\uCA^1_{q1}
\,\bulpsi^{b\kern 0.4pt\alpha}_{111}[i],
\qquad
\mytag{6.5}\\
\vspace{2ex}
&\nabla_{\!q}\,\circpsi^{\,b}_{111111}[i]
\to\nabla_{\!q}\,\circpsi^{\,b}_{111111}[i]
+\frac{6\,i\,e\,g_1}{\hbar\,c}\,\uCA^1_{q1}
\,\circpsi^{\,b}_{111111}[i].
\mytag{6.6}
\endalign
$$
Before substituting \mythetag{6.5} and \mythetag{6.6} back into
\mythetag{6.1} let's apply \mythetag{4.35} and \mythetag{6.3} to
\mythetag{6.5}. As a result we get the following expression:
$$
\allowdisplaybreaks
\gather
\nabla_{\!q}\bulpsi^{b\kern 0.4pt\alpha}_{111}[i]
\to\nabla_{\!q}\bulpsi^{\,b}_{111111}[i]
\cdot\frac{\varphi^{\kern 0.4pt\alpha111}[vac]}{|\boldsymbol
\varphi[vac]|}+\nabla_{\!q}\bulpsi^{\,b}[i]\cdot\frac{\phi^{\kern
0.8pt\alpha}_{111}[vac]}{|\boldsymbol\phi[vac]|}\,-\\
\vspace{2ex}
-\,\frac{i\,e\,g_2}{\hbar\,c}\left(\!\frac{3\,g_1\,A_q-g_2\,Z_q}
{\sqrt{(g_2)^2+(3\,g_1)^2}}\left(\bulpsi^{\,b}[i]\cdot\frac{\phi^{\kern
0.8pt\alpha}_{111}[vac]}{|\boldsymbol\phi[vac]|}
-\bulpsi^{\,b}_{111111}[i]\,\cdot\right.\right.\\
\vspace{2ex}
\left.\left.\cdot\,\frac{\varphi^{\kern 0.4pt\alpha111}[vac]}
{|\boldsymbol\varphi[vac]|}
\right)+W^{{\sssize +}111111}_q\,\bulpsi^{\,b}_{111111}[i]
\cdot\frac{\phi^{\kern 0.8pt\alpha}_{111}[vac]}
{|\boldsymbol\phi[vac]|}+W^{\sssize -}_{q111111}\,\times\right.\\
\vspace{2ex}
\left.\times\,\bulpsi^{\,b}[i]\cdot\frac{\varphi^{\kern 0.4pt
\alpha111}[vac]}{|\boldsymbol\varphi[vac]|}
\vphantom{\frac{3\,g_1\,A_q-g_2\,Z_q}
{\sqrt{(g_2)^2+(3\,g_1)^2}}}\right)+\frac{3\,i\,e\,g_1}{\hbar\,c}
\,\frac{g_2\,A_q+3\,g_1\,Z_q}{\sqrt{(g_2)^2+(3\,g_1)^2}}\,\times\\
\vspace{2ex}
\times\left(\bulpsi^{\,b}_{111111}[i]
\cdot\frac{\varphi^{\kern 0.4pt\alpha111}[vac]}{|\boldsymbol
\varphi[vac]|}+
\bulpsi^{\,b}[i]\cdot\frac{\phi^{\kern 0.8pt\alpha}_{111}[vac]}
{|\boldsymbol\phi[vac]|}
\right).
\endgather
$$
In deriving the above expression for $\nabla_{\!q}\bulpsi^{b\kern
0.4pt\alpha}_{111}[i]$ we used the formula \mythetag{4.33} for
$\uCA^1_{k1}$ and the formula \mythetag{4.35} for $\uuCA^\alpha_{q
\kern 0.4pt\theta}$. Recollecting terms in it, we get
$$
\gathered
\nabla_{\!q}\bulpsi^{b\kern 0.4pt\alpha}_{111}[i]
\to\left(\nabla_{\!q}\bulpsi^{\,b}_{111111}[i]
+\frac{i\,e\,g_2}{\hbar\,c}\,\frac{3\,g_1\,A_q-g_2\,Z_q}
{\sqrt{(g_2)^2+(3\,g_1)^2}}\ \bulpsi^{\,b}_{111111}[i]\,+\right.\\
\vspace{2ex}
\left.+\frac{3\,i\,e\,g_1}{\hbar\,c}\,\frac{g_2\,A_q+3\,g_1\,Z_q}
{\sqrt{(g_2)^2+(3\,g_1)^2}}\ \bulpsi^{\,b}_{111111}[i]
-\frac{i\,e\,g_2}{\hbar\,c}\,W^{\sssize -}_{q111111}\,\bulpsi^{\,b}[i]
\right)\cdot\\
\vspace{2ex}
\cdot\,\frac{\varphi^{\kern 0.4pt\alpha111}[vac]}{|\boldsymbol
\varphi[vac]|}+\left(\nabla_{\!q}\bulpsi^{\,b}[i]
-\frac{i\,e\,g_2}{\hbar\,c}\,\frac{3\,g_1\,A_q-g_2\,Z_q}
{\sqrt{(g_2)^2+(3\,g_1)^2}}\ \bulpsi^{\,b}[i]
+\frac{3\,i\,e\,g_1}{\hbar\,c}\,\times\right.\\
\vspace{2ex}
\left.\times\,\frac{g_2\,A_q+3\,g_1\,Z_q}{\sqrt{(g_2)^2+(3\,g_1)^2}}
\,\bulpsi^{\,b}[i]-\frac{i\,e\,g_2}{\hbar\,c}
\,W^{{\sssize +}111111}_q\,\bulpsi^{\,b}_{111111}[i]\right)
\cdot\frac{\phi^{\kern 0.8pt\alpha}_{111}[vac]}{|\boldsymbol
\phi[vac]|}.
\endgathered\quad
\mytag{6.7}
$$
Simplifying the right hand side of the formula \mythetag{6.7} 
a little bit more, we find
$$
\gathered
\nabla_{\!q}\bulpsi^{b\kern 0.4pt\alpha}_{111}[i]
\to\left(\nabla_{\!q}\bulpsi^{\,b}_{111111}[i]
+\frac{i\,e}{\hbar\,c}\,\frac{6\,g_1\,g_2\,A_q}
{\sqrt{(g_2)^2+(3\,g_1)^2}}\ \bulpsi^{\,b}_{111111}[i]\,+\right.\\
\vspace{2ex}
\left.+\frac{i\,e}{\hbar\,c}\,\frac{(3\,g_1)^2-(g_2)^2}
{\sqrt{(g_2)^2+(3\,g_1)^2}}\ Z_q\,\bulpsi^{\,b}_{111111}[i]
-\frac{i\,e\,g_2}{\hbar\,c}\,W^{\sssize -}_{q111111}\,\bulpsi^{\,b}[i]
\right)\cdot\\
\vspace{2ex}
\cdot\,\frac{\varphi^{\kern 0.4pt\alpha111}[vac]}{|\boldsymbol
\varphi[vac]|}+\left(\nabla_{\!q}\bulpsi^{\,b}[i]
+\frac{i\,e}{\hbar\,c}\,\sqrt{(g_2)^2+(3\,g_1)^2}\,Z_q\, 
\bulpsi^{\,b}[i]\,-\right.\\
\vspace{2ex}
\left.-\,\frac{i\,e\,g_2}{\hbar\,c}
\,W^{{\sssize +}111111}_q\,\bulpsi^{\,b}_{111111}[i]\right)
\cdot\frac{\phi^{\kern 0.8pt\alpha}_{111}[vac]}{|\boldsymbol
\phi[vac]|}.
\endgathered\quad
\mytag{6.8}
$$
Acting in a similar way, from \mythetag{6.6} we derive the following
formula:
$$
\pagebreak
\gathered
\nabla_{\!q}\,\circpsi^{\,b}_{111111}[i]
\to\nabla_{\!q}\,\circpsi^{\,b}_{111111}[i]
+\frac{i\,e}{\hbar\,c}\,\frac{6\,g_1\,g_2\,A_q}
{\sqrt{(g_2)^2+(3\,g_1)^2}}\,\times\\
\vspace{2ex}
\times\,\circpsi^{\,b}_{111111}[i]+\frac{i\,e}{\hbar\,c}\,
\frac{18\,(g_1)^2}{\sqrt{(g_2)^2+(3\,g_1)^2}}\ Z_q
\,\circpsi^{\,b}_{111111}[i].
\endgathered
\mytag{6.9}
$$\par
     The expansion \mythetag{6.4} can be teated as the expansion of
lepton wave functions into chiral and antichiral parts. Therefore,
using \mythetag{1.10}, we write
$$
\xalignat 2
&\hskip -2em
\bulpsi^{\,b}_{111111}[i]=\sum^2_{c=1}\bulletH^b_c\ \psi^c_{111111}[i],
&&\circpsi^{\,b}_{111111}[i]=\sum^2_{c=1}\circH^b_c\ \psi^c_{111111}[i].
\qquad
\mytag{6.10}
\endxalignat
$$
Here $\bulletH^b_c$ and $\circH^b_c$ are the components of the chiral 
and antichiral projection operators introduced in \mythetag{1.9} and
\mythetag{1.10}.\par
     Now let's substitute \mythetag{6.8} and \mythetag{6.9} into
\mythetag{6.1}. As a result, using the formulas \mythetag{2.64},
\mythetag{4.6}, \mythetag{4.7} and the orthogonality condition
\mythetag{4.10}, we derive
$$
\hskip -2em
\Cal L_6=\Cal L_{61}+\Cal L_{62}+\Cal L_{63}
+\Cal L_{64}+\Cal L_{65}.
\mytag{6.11}
$$
The first term $\Cal L_{61}$ in the right hand side of the formula
\mythetag{6.11} is a standard kinetic term for three spin $1/2$ 
particles with electric charge $Q$ in an electromagnetic field:
$$
\aligned
\Cal L_{61}&=i\,\hbar\int\!\!\sum_{i=e,\mu,\tau}
\sum^4_{a=1}\sum^4_{\bar a=1}\sum^4_{b=1}\sum^3_{q=0}
\uD^{11}\,\uD^{11}\,\uD^{11}\,\uD^{11}\,\uD^{11}\,\uD^{11}
\,\times\\
\vspace{2ex}
&\times\,
\overline{\psi^{\raise 0.6pt \hbox{$\ssize\bar a$}}_{111111}[i]}
\,D_{a\bar a}\,\gamma^{\kern 0.5pt aq}_b\left(\!\nabla_{\!q}
\,\psi^{\,b}_{111111}[i]-\frac{i\,Q}{\hbar\,c}\,A_q
\,\psi^{\,b}_{111111}[i]\right)\,dV.
\endaligned\quad
\mytag{6.12}
$$ 
The electric charge $Q$ for all charged leptons is given by the 
following formula:
$$
\hskip -2em
Q=-\frac{6\,e\,g_1\,g_2}{\sqrt{(g_2)^2+(3\,g_1)^2}}.
\mytag{6.13}
$$
Since one of the three particles described by \mythetag{6.12} is an
electron, its charge $Q=-e$. Therefore, from \mythetag{6.13} we derive 
the following equality relating $g_1$ and $g_2$:
$$
\hskip -2em
\frac{6\,g_1\,g_2}{\sqrt{(g_2)^2+(3\,g_1)^2}}=1.
\mytag{6.14}
$$
Due to \mythetag{6.13} and \mythetag{6.14} the electric charge of all 
three charged leptons $e$, $\mu$, and $\tau$ in the table \mythetag{1.1} 
is negative and is equal to the charge of an electron:
$$
\hskip -2em
Q=-e.
\mytag{6.15}
$$\par
    The term $\Cal L_{62}$ in \mythetag{6.11} is a standard kinetic 
term describing three electrically neutral leptons --- $e$-neutrino, 
$\mu$-neutrino, abd $\tau$-neutrino:
$$
\Cal L_{62}=i\,\hbar\int\!\!\sum_{i=e,\mu,\tau}
\sum^4_{a=1}\sum^4_{\bar a=1}\sum^4_{b=1}\sum^3_{q=0}
\overline{\bulpsi^{\raise 0.6pt \hbox{$\ssize\bar a$}}[i]}
\,D_{a\bar a}\,\gamma^{\kern 0.5pt aq}_b\,\nabla_{\!q}
\bulpsi^{\,b}[i]\,dV.\quad
\mytag{6.16}
$$ 
Note that the wave functions of neutrinos \mythetag{6.16} have only 
chiral components. Their antichiral components are zero. This is the 
sign of chiral-to-antichiral asymmetry \pagebreak of the Standard 
Model.\par
     The term $\Cal L_{63}$ in \mythetag{6.11} is a purely potential
term. It describes the interaction of charged leptons with $Z$-bosons.
Applying \mythetag{6.10}, we write
$$
\gathered
\Cal L_{63}=-\frac{e}{c}\int\!\!\sum_{i=e,\mu,\tau}
\sum^4_{a=1}\sum^4_{\bar a=1}\sum^4_{b=1}\sum^4_{c=1}\sum^3_{q=0}
\uD^{11}\,\uD^{11}\,\uD^{11}\,\uD^{11}\,\uD^{11}\,\uD^{11}\,
\times\\
\vspace{2ex}
\times\,\overline{\psi^{\raise 0.6pt \hbox{$\ssize\bar a$}}_{111111}[i]}
\,D_{a\bar a}\,\gamma^{\kern 0.5pt aq}_c\,Z_q
\frac{((3\,g_1)^2-(g_2)^2)\,\bulletH^c_b+18\,(g_1)^2\,\circH^c_b}
{\sqrt{(g_2)^2+(3\,g_1)^2}}\ \psi^{\,b}_{111111}[i]\,dV.
\endgathered\quad
\mytag{6.17}
$$
The presence of the components of chiral and anichiral projectors
\mythetag{1.9} in \mythetag{6.17} is another sign of chiral-to-antichiral
asymmetry of the Standard Model.\par
     The term $\Cal L_{64}$ in \mythetag{6.11} is a purely potential
term. It describes the interaction of electrically neutral leptons 
with $Z$-bosons. From \mythetag{6.8} we derive
$$
\gathered
\Cal L_{64}=-\frac{e}{c}\int\!\!\sum_{i=e,\mu,\tau}
\sum^4_{a=1}\sum^4_{\bar a=1}\sum^4_{b=1}\sum^4_{c=1}\sum^3_{q=0}
\overline{\psi^{\raise 0.6pt \hbox{$\ssize\bar a$}}[i]}
\,D_{a\bar a}\,\times\\
\vspace{2ex}
\times\,\gamma^{\kern 0.5pt aq}_c\,Z_q
\sqrt{(g_2)^2+(3\,g_1)^2}\,\bulletH^c_b
\ \psi^{\,b}[i]\,dV.
\endgathered\quad
\mytag{6.18}
$$
Like \mythetag{6.17}, this term \mythetag{6.18} also breaks the
chiral-to-antichiral symmetry.\par
    The term $\Cal L_{65}$ in the sum \mythetag{6.11} is also a purely
potential term. It describes the interaction of charged and neutral 
leptons with $W$-bosons:
$$
\gathered
\Cal L_{65}=\frac{e}{c}\int\!\!\sum_{i=e,\mu,\tau}
\sum^4_{a=1}\sum^4_{\bar a=1}\sum^4_{b=1}\sum^3_{q=0}
\uD^{11}\,\uD^{11}\,\uD^{11}\,\uD^{11}\,\uD^{11}\,\uD^{11}\,
\times\\
\vspace{2ex}
\times\,\overline{\psi^{\raise 0.6pt \hbox{$\ssize\bar a$}}_{111111}[i]}
\,D_{a\bar a}\,\gamma^{\kern 0.5pt aq}_b\,g_2
\,W^{\sssize -}_{q111111}\,\bulpsi^{\,b}[i]\ dV\,+\\
\vspace{2ex}
+\,\frac{e}{c}\int\!\!\sum_{i=e,\mu,\tau}
\sum^4_{a=1}\sum^4_{\bar a=1}\sum^4_{b=1}\sum^3_{q=0}
\overline{\bulpsi^{\raise 0.6pt \hbox{$\ssize\bar a$}}[i]}
\,D_{a\bar a}\,\gamma^{\kern 0.5pt aq}_b\,g_2
\,W^{{\sssize +}111111}_q\,\psi^{\,b}_{111111}[i]\ dV.
\endgathered\quad
\mytag{6.19}
$$
The integral $\Cal L_{65}$ in \mythetag{6.19} has the same sign of
chiral-to-antichiral asymmetry of the Standard Model as the integral
\mythetag{6.16}.\par
     Now let's proceed with the integral \mythetag{6.2}. It is a 
purely potential action integral. Substituting \mythetag{6.3} into
\mythetag{6.2}, we take into account \mythetag{2.64}, \mythetag{4.6}, 
\mythetag{4.7}, and the orthogonality condition \mythetag{4.10}. 
Apart from those mentioned above, we take into account the formula
\mythetag{3.12}. As a result we get
$$
\hskip -2em
\Cal L_7=\Cal L_{71}+\Cal L_{72}.
\mytag{6.20}
$$
The first term $\Cal L_{71}$ in the expansion \mythetag{6.20} is a purely
potential term. Moreover, it is a mass term. \pagebreak It determines the
masses of three charged leptons --- an electron, a muon, and a tauon in 
the leptons generation table \mythetag{1.1}:
$$
\gathered
\Cal L_{71}=-\!\!\sum_{i=e,\mu,\tau}\!\!\frac{h[i]\,v}{\sqrt{2}}
\!\int\!\sum^4_{a=1}
\sum^4_{\bar a=1}\sum^2_{\alpha=1}\sum^2_{\bar\alpha=1}
\uD^{11}\,\uD^{11}\,\uD^{11}\,\uD^{11}\,\uD^{11}\,\uD^{11}\,\times\\
\vspace{2ex}
\times\,D_{a\bar a}\ \overline{\psi^{\raise 0.6pt \hbox{$\ssize
\bar a$}}_{111111}[i]}\ \psi^a_{111111}[i]\ dV.
\endgathered
\quad
\mytag{6.21}
$$
The second term $\Cal L_{72}$ in \mythetag{6.20} is very similar to
\mythetag{6.21}. It is also a purely potential term describing the
interaction of charged leptons with the real scalar Higgs field:
$$
\gathered
\Cal L_{72}=-\!\sum_{i=e,\mu,\tau}\!\frac{h[i]}{\sqrt{2}}\int\!
\sum^4_{a=1}
\sum^4_{\bar a=1}\sum^2_{\alpha=1}\sum^2_{\bar\alpha=1}
\uD^{11}\,\uD^{11}\,\uD^{11}\,\uD^{11}\,\uD^{11}\,\uD^{11}\,\times\\
\vspace{2ex}
\times\,D_{a\bar a}\ \chi\ \overline{\psi^{\raise 0.6pt \hbox{$\ssize
\bar a$}}_{111111}[i]}\ \psi^a_{111111}[i]\ dV.
\endgathered
\quad
\mytag{6.22}
$$
From \mythetag{6.21} we derive the following formulas for the masses
of charged leptons:
$$
\xalignat 3
&\hskip -2em
m_e=\frac{h[e]\,v}{\sqrt{2}\,c}\,,
&&m_\mu=\frac{h[\mu]\,v}{\sqrt{2}\,c}\,,
&&m_\tau=\frac{h[\tau]\,v}{\sqrt{2}\,c}\,.
\qquad
\mytag{6.23}
\endxalignat
$$
The lepton part of the total action integral is exhausted by
\mythetag{6.21} and \mythetag{6.22}. For this reason uncharged
leptons are massless particles in the Standard Model.
\head
7. The quark masses.
\endhead
     Quarks are represented in the table \mythetag{1.2}. Like
leptons they are subdivided into three pairs (three generations).
Like in the case of leptons, the chiral parts of quark wave
functions form $\MatGrSU(2)$-doublets, while their antichiral
parts are singlets. Here is the kinetic term of the quark action
integral:
$$
\hskip -2em
\aligned
\Cal L_8&=i\,\hbar\int\sum^3_{i=1}
\sum^4_{a=1}\sum^4_{\bar a=1}\sum^4_{b=1}\sum^3_{q=0}
\sum^2_{\alpha=1}\sum^2_{\bar\alpha=1}
\sum^3_{\beta=1}\sum^3_{\bar\beta=1}
\uD_{11}\,\times\\
\vspace{2ex}
&\qquad\times\,\uuD_{\alpha\bar\alpha}
\,\uuuD_{\beta\bar\beta}\,
\overline{\bulpsi^{\raise 0.6pt \hbox{$\ssize\bar a1
\bar\alpha\bar\beta$}}[i]}\,D_{a\bar a}\,\gamma^{\kern 0.5pt aq}_b
\ \nabla_{\!q}\bulpsi^{\kern 0.5pt b\kern 0.2pt 1\alpha\beta}[i]
\ dV\,+\\
\vspace{2ex}
&+\,i\,\hbar\int\!\!\sum_{i=u,c,t}
\sum^4_{a=1}\sum^4_{\bar a=1}\sum^4_{b=1}\sum^3_{q=0}
\sum^2_{\alpha=1}\sum^2_{\bar\alpha=1}
\sum^3_{\beta=1}\sum^3_{\bar\beta=1}
\uD_{11}\,\uD_{11}\,\uD_{11}\,\times\\
\vspace{2ex}
&\times\,\uD_{11}\,
\uuD_{\alpha\bar\alpha}
\,\uuuD_{\beta\bar\beta}\,
\overline{\circpsi^{\raise 0.6pt \hbox{$\ssize\bar a1111\bar\beta$}}[i]}
\,D_{a\bar a}\,\gamma^{\kern 0.5pt aq}_b
\ \nabla_{\!q}\circpsi^{b\kern 0.2pt 1111\beta}[i]\,dV\,+\\
\vspace{2ex}
&+\,i\,\hbar\int\!\!\sum_{i=d,s,b}
\sum^4_{a=1}\sum^4_{\bar a=1}\sum^4_{b=1}\sum^3_{q=0}
\sum^2_{\alpha=1}\sum^2_{\bar\alpha=1}
\sum^3_{\beta=1}\sum^3_{\bar\beta=1}
\uD^{11}\,\times\\
\vspace{2ex}
&\qquad\times\,\uD^{11}\,\uuD_{\alpha\bar\alpha}
\,\uuuD_{\beta\bar\beta}\,
\overline{\circpsi^{\raise 0.6pt \hbox{$\ssize\bar a\bar\beta$}}_{11}[i]}
\,D_{a\bar a}\,\gamma^{\kern 0.5pt aq}_b
\ \nabla_{\!q}\circpsi^{b\beta}_{11}[i]\,dV.
\endaligned
\mytag{7.1}
$$ 
The mass terms for the quark action integral are more complicated 
as compared to the case of leptons because of the generation mixing:
$$
\hskip -2em
\aligned
\Cal L_9&=-\sum^3_{i=1}\sum^3_{j=1}h_1[ij]\!\int\!\sum^4_{a=1}
\sum^4_{\bar a=1}\sum^2_{\alpha=1}\sum^2_{\bar\alpha=1}
\sum^3_{\beta=1}\sum^3_{\bar\beta=1}
\uD_{11}\,\times\\
\vspace{2ex}
&\qquad\times\,\uuD_{\alpha\bar\alpha}
\,\uuuD_{\beta\bar\beta}\,D_{a\bar a}
\,\overline{\circpsi^{\raise 0.6pt \hbox{$\ssize\bar a1111
\bar\beta$}}[i]}\ \overline{\phi^{\kern 0.4pt
\raise 1.2pt\hbox{$\ssize\kern 0.8pt\bar\alpha$}}_{111}}
\,\bulpsi^{\kern 0.5pt a1\alpha\beta}[j]\,dV\,-\\
\vspace{2ex}
&-\sum^3_{i=1}\sum^3_{j=1}\overline{h_1[j\kern 0.5pt i]
\vphantom{A^A}}\!\int\!\sum^4_{a=1}\sum^4_{\bar a=1}
\sum^2_{\alpha=1}\sum^2_{\bar\alpha=1}\sum^3_{\beta=1}
\sum^3_{\bar\beta=1}\uD_{11}\,\times\\
\vspace{2ex}
&\qquad\times\,\uuD_{\alpha\bar\alpha}
\,\uuuD_{\beta\bar\beta}\,D_{a\bar a}
\,\overline{\bulpsi^{\raise 0.6pt \hbox{$\kern 0.5pt\ssize\bar a1
\bar\alpha\bar\beta$}}[i]}\ \phi^{\kern 1.2pt\alpha}_{111}
\,\circpsi^{a1111\beta}[j]\,dV\,-\\
\vspace{2ex}
&-\sum^3_{i=1}\sum^3_{j=1}h_2[ij]\!\int\!\sum^4_{a=1}
\sum^4_{\bar a=1}\sum^2_{\alpha=1}\sum^2_{\bar\alpha=1}
\sum^3_{\beta=1}\sum^3_{\bar\beta=1}\uD_{11}\,\times\\
\vspace{2ex}
&\qquad\qquad\times\,\uuD_{\alpha\bar\alpha}
\,\uuuD_{\beta\bar\beta}\,D_{a\bar a}
\,\overline{\circpsi^{\raise 0.6pt \hbox{$\kern 0.5pt\ssize
\bar a\bar\beta$}}_{11}[i]}\ \overline{\varphi^{\raise 
1.2pt\hbox{$\ssize\bar\alpha111$}}}
\,\bulpsi^{\kern 0.5pt a1\alpha\beta}[j]\,dV\,-\\
\vspace{2ex}
&-\sum^3_{i=1}\sum^3_{j=1}\overline{h_2[j\kern 0.5pt i]
\vphantom{A^A}}\!\int\!\sum^4_{a=1}
\sum^4_{\bar a=1}\sum^2_{\alpha=1}\sum^2_{\bar\alpha=1}
\sum^3_{\beta=1}\sum^3_{\bar\beta=1}\uD_{11}\,\times\\
\vspace{2ex}
&\qquad\qquad\times\,\uuD_{\alpha\bar\alpha}
\,\uuuD_{\beta\bar\beta}\,D_{a\bar a}
\,\overline{\bulpsi^{\raise 0.6pt \hbox{$\kern 0.5pt\ssize
\bar a1\bar\alpha\bar\beta$}}[i]}
\ \varphi^{\alpha111}\,\circpsi^{\kern 0.5pt a\beta}_{11}[j]
\,dV.
\endaligned
\mytag{7.2}
$$
Like in \mythetag{6.1} and in \mythetag{6.2} for leptons, in the
above two action integrals for quarks both singlet and doublet wave
functions \mythetag{1.13} and \mythetag{1.14} are used. But instead
of three real parameters $h[e]$, $h[\mu]$, $h[\tau]$ here we have
two complex $3\times 3$ matrices of such parameters $h_1[ij]$ and
$h_2[ij]$. In the case of leptons we used the expansion 
\mythetag{6.3}. For quarks such an expansion looks like
$$
\hskip -2em
\bulpsi^{\kern 0.5pt a\kern 0.2pt 1\alpha\beta}[i]
=\bulpsi^{\kern 0.5pt a\beta}_{11}[i]
\cdot\frac{\varphi^{\kern 0.4pt\alpha111}[vac]}{|\boldsymbol
\varphi[vac]|}+\bulpsi^{a1111\beta}[i]\cdot\frac{\phi^{\kern
0.8pt\alpha}_{111}[vac]}{|\boldsymbol\phi[vac]|}.
\mytag{7.3}
$$
The chiral coefficients $\bulpsi^{\kern 0.5pt a\beta}_{11}[i]$ and
$\bulpsi^{a1111\beta}[i]$ from \mythetag{7.3} are complementary to
the antichiral wave functions \mythetag{1.13} and \mythetag{1.14}.
Indeed, we can define
$$
\hskip -2em
\aligned
&\psi^{a1111\beta}[i]=\bulpsi^{a1111\beta}[i]+\circpsi^{a1111\beta}[i],\\
\vspace{2ex}
&\psi^{\kern 0.5pt a\beta}_{11}[i]=\bulpsi^{\kern 0.5pt a\beta}_{11}[i]
+\circpsi^{\kern 0.5pt a\beta}_{11}[i].
\endaligned
\mytag{7.4}
$$
The formulas \mythetag{7.4} are analogous to \mythetag{6.4}. \pagebreak 
The wave functions introduced by means of the formulas \mythetag{7.4} 
are interpreted as the complete wave functions of quarks:
$$
\xalignat 3
&\hskip -2em
\psi^{a1111\beta}[u],
&&\psi^{a1111\beta}[c],
&&\psi^{a1111\beta}[t],\\
\vspace{-1ex}
&&&\mytag{7.5}\\
\vspace{-1ex}
&\hskip -2em
\psi^{a\beta}_{11}[d],
&&\psi^{a\beta}_{11}[s],
&&\psi^{a\beta}_{11}[b].
\endxalignat
$$
The first line in \mythetag{7.5} corresponds to upper level quarks,
i\.\,e\. an up-quark, a charm-quark, and a top-quark, the second line
describes lower level quarks --- a down-quark, a strange-quark, and a 
bottom-quark.\par
     The covariant derivatives $\nabla_q$ in \mythetag{7.1} are complete
non-vacuum covariant derivatives. They are evaluated according to the
formulas \mythetag{2.5}, \mythetag{2.6}, and \mythetag{2.7}, where 
the electro-weak connection components are taken from \mythetag{2.36},
\mythetag{2.37}, and \mythetag{2.38} including the components of the 
gauge fields $\bolduCA$,  $\bolduuCA$, and $\bolduuuCA$. Passing from
non-vacuum to vacuum covariant derivatives in \mythetag{2.5},
\mythetag{2.6}, and \mythetag{2.7} we get
$$
\align
&\hskip -2em
\aligned
\nabla_{\!q}&\bulpsi^{a1\alpha\beta}[i]
\to\nabla_{\!q}\bulpsi^{a1\alpha\beta}[i]
-\frac{i\,e\,g_1}{\hbar\,c}\,\uCA^1_{q1}
\,\bulpsi^{a1\alpha\beta}[i]\,-\\
&-\,\frac{i\,e\,g_2}{\hbar\,c}\sum^2_{\theta=1}
\uuCA^\alpha_{q\kern 0.4pt\theta}\,
\bulpsi^{a1\theta\kern 0.4pt\beta}[i]
-\frac{i\,e\,g_3}{\hbar\,c}\sum^3_{\theta=1}
\uuuCA^\alpha_{q\kern 0.4pt\theta}
\,\bulpsi^{a1\alpha\kern 0.4pt\theta}[i],
\endaligned
\mytag{7.6}\\
\vspace{2ex}
&\hskip -2em
\aligned
\nabla_{\!q}&\circpsi^{a1111\beta}[i]
\to\nabla_{\!q}\circpsi^{a1111\beta}[i]\,-\\
&-\,\frac{4\,i\,e\,g_1}{\hbar\,c}\,\uCA^1_{q1}
\,\circpsi^{a1111\beta}[i]
-\frac{i\,e\,g_3}{\hbar\,c}\sum^3_{\theta=1}
\uuuCA^\alpha_{q\kern 0.4pt\theta}
\,\circpsi^{a1111\kern 0.4pt\theta}[i],
\endaligned
\mytag{7.7}\\
\vspace{2ex}
&\hskip -2em
\aligned
\nabla_{\!q}&\circpsi^{a\beta}_{11}[i]
\to\nabla_{\!q}\circpsi^{a\beta}_{11}[i]
+\frac{2\,i\,e\,g_1}{\hbar\,c}\,\uCA^1_{q1}
\,\circpsi^{a\beta}_{11}[i]\,-\\
&-\,\frac{i\,e\,g_3}{\hbar\,c}\sum^3_{\theta=1}\uuuA^\alpha_{q
\kern 0.4pt\theta}\,\circpsi^{a\kern 0.4pt\theta}_{11}[i].
\endaligned
\mytag{7.8}
\endalign
$$
The next step is to apply the formulas \mythetag{4.33} and \mythetag{4.35}
to \mythetag{7.6}, \mythetag{7.7}, and \mythetag{7.8}. Let's begin with 
the last two formulas \mythetag{7.7} and \mythetag{7.8} for 
$\MatGrSU(2)$-singlet functions. In the case of the formula \mythetag{7.7},
applying \mythetag{4.33} to it, we derive
$$
\aligned
\nabla_{\!q}&\circpsi^{a1111\beta}[i]
\to\nabla_{\!q}\circpsi^{a1111\beta}[i]
-\frac{i\,e\,}{\hbar\,c}\,\frac{4\,g_1\,g_2\,A_q}
{\sqrt{(g_2)^2+(3\,g_1)^2}}\ \circpsi^{a1111\beta}[i]\,-\\
\vspace{2ex}
&-\,\frac{i\,e}{\hbar\,c}\,\frac{12\,(g_1)^2\,Z_q}
{\sqrt{(g_2)^2+(3\,g_1)^2}}\ \circpsi^{a1111\beta}[i]
-\frac{i\,e\,g_3}{\hbar\,c}\sum^3_{\theta=1}
\uuuCA^\alpha_{q\kern 0.4pt\theta}
\,\circpsi^{a1111\kern 0.4pt\theta}[i].
\endaligned\quad
\mytag{7.9}
$$
In a similar way, applying \mythetag{4.33} to the formula
\mythetag{7.8}, we derive
$$
\hskip -2em
\aligned
\nabla_{\!q}&\circpsi^{a\beta}_{11}[i]
\to\nabla_{\!q}\circpsi^{a\beta}_{11}[i]
+\frac{i\,e}{\hbar\,c}\,
\frac{2\,g_1\,g_2\,A_k}{\sqrt{(g_2)^2+(3\,g_1)^2}}
\ \circpsi^{a\beta}_{11}[i]\,+\\
\vspace{2ex}
&+\,\frac{i\,e}{\hbar\,c}\,
\frac{6\,(g_1)^2}{\sqrt{(g_2)^2+(3\,g_1)^2}}\ Z_q
\,\circpsi^{a\beta}_{11}[i]
-\frac{i\,e\,g_3}{\hbar\,c}\sum^3_{\theta=1}\uuuA^\alpha_{q
\kern 0.4pt\theta}\,\circpsi^{a\kern 0.4pt\theta}_{11}[i].
\endaligned\quad
\mytag{7.10}
$$
Applying \mythetag{4.33} and \mythetag{4.35} to \mythetag{7.6}, we 
should take into account the formula \mythetag{7.3}:
$$
\gathered
\nabla_{\!q}\bulpsi^{a1\alpha\beta}[i]
\to\left(\nabla_{\!q}\bulpsi^{\kern 0.5pt a\beta}_{11}[i]
+\frac{i\,e}{\hbar\,c}\,\frac{2\,g_1\,g_2\,A_q}{\sqrt{(g_2)^2
+(3\,g_1)^2}}\ \bulpsi^{\kern 0.5pt a\beta}_{11}[i]
\,-\right.\\
\vspace{2ex}
-\,\frac{i\,e}{\hbar\,c}\,\frac{(3\,(g_1)^2+(g_2)^2)}
{\sqrt{(g_2)^2+(3\,g_1)^2}}\ Z_q\,\bulpsi^{\kern 0.5pt
a\beta}_{11}[i]-\,\frac{i\,e\,g_2}
{\hbar\,c}\,W^{\sssize -}_{q111111}\,\bulpsi^{a1111\beta}[i]\,-\\
\vspace{2ex}
\left.-\,\frac{i\,e\,g_3}{\hbar\,c}\,\sum^3_{\theta=1}
\uuuCA^\alpha_{q\kern 0.4pt\theta}\,\bulpsi^{\kern 0.5pt a
\theta}_{11}[i]\right)\cdot\,\frac{\varphi^{\kern 0.4pt
\alpha111}[vac]}{|\boldsymbol\varphi[vac]|}+\left(\nabla_{\!q}
\bulpsi^{a1111\beta}[i]-\frac{i\,e}{\hbar\,c}\,\times
\vphantom{\frac{3\,g_1\,A_q-g_2\,Z_q}
{\sqrt{(g_2)^2+(3\,g_1)^2}}}\right.\\
\vspace{2ex}
\times\,\frac{4\,g_1\,g_2\,A_q}
{\sqrt{(g_2)^2+(3\,g_1)^2}}\,\bulpsi^{a1111\beta}[i]
-\frac{i\,e}{\hbar\,c}\,\frac{3\,(g_1)^2-(g_2)^2}
{\sqrt{(g_2)^2+(3\,g_1)^2}}\,Z_q\,\bulpsi^{a1111\beta}[i]\,-\\
\vspace{2ex}
\left.-\,\frac{i\,e\,g_2}{\hbar\,c}
\,W^{{\sssize +}111111}_q\,\bulpsi^{\kern 0.5pt a\beta}_{11}[i]
-\frac{i\,e\,g_3}{\hbar\,c}\,\sum^3_{\theta=1}
\uuuCA^\alpha_{q\kern 0.4pt\theta}\,\bulpsi^{a1111\theta}[i]
\right)\cdot\frac{\phi^{\kern 0.8pt\alpha}_{111}[vac]}
{|\boldsymbol\phi[vac]|}.\\
\endgathered\quad
\mytag{7.11}
$$
The formula \mythetag{7.11} is analogous to \mythetag{6.8}. When
substituting \mythetag{7.9}, \mythetag{7.10}, and \mythetag{7.11}
into \mythetag{7.1} we take into account the equality \mythetag{6.14}.
Then we get 
$$
\hskip -2em
\Cal L_8=\Cal L_{81}+\Cal L_{82}+\Cal L_{83}+\Cal L_{84}+\Cal L_{85}.
\mytag{7.12}
$$
The first term $\Cal L_{81}$ in the right hand side of the formula
\mythetag{7.12} is a standard kinetic term for three spin $1/2$ 
particles in an electromagnetic field:
$$
\aligned
\Cal L_{81}&=i\,\hbar\int\!\!\sum_{i=d,s,b}
\ \sum^4_{a=1}\sum^4_{\bar a=1}\sum^4_{b=1}\sum^3_{q=0}
\sum^3_{\beta=1}\sum^3_{\bar\beta=1}
\uD^{11}\,\uD^{11}\,\uuuD_{\beta\bar\beta}
\ \overline{\psi^{\raise 0.6pt \hbox{$\ssize\kern 0.5pt
\bar a\bar\beta$}}_{11}[i]}\,D_{a\bar a}
\,\gamma^{\kern 0.5pt aq}_b\,\times\\
\vspace{2ex}
&\times\left(\!
\nabla_{\!q}\psi^{\kern 0.5pt b\beta}_{11}[i]
+\frac{i\,e}{3\,\hbar\,c}\,A_q\,
\psi^{\kern 0.5pt b\beta}_{11}[i]
-\frac{i\,e\,g_3}{\hbar\,c}\,\sum^3_{\theta=1}
\uuuCA^\alpha_{q\kern 0.4pt\theta}\,\psi^{\kern 0.5pt
b\kern 0.2pt\theta}_{11}[i]\right)\,dV.
\endaligned\quad
\mytag{7.13}
$$ 
The second term $\Cal L_{82}$ is a standard kinetic term for
other three particles
$$
\hskip -2em
\aligned
\Cal L_{82}&=i\,\hbar\int\!\!\sum_{i=u,c,t}\,
\sum^4_{a=1}\sum^4_{\bar a=1}\sum^4_{b=1}\sum^3_{q=0}
\sum^3_{\beta=1}\sum^3_{\bar\beta=1}
\uD^{11}\,\uD^{11}\,\times\\
\vspace{2ex}
&\times\,\uD^{11}\,\uD^{11}\,\uuuD_{\beta\bar\beta}
\ \overline{\psi^{\raise 0.6pt \hbox{$\ssize\kern 0.5pt
\bar a1111\bar\beta$}}[i]}\,D_{a\bar a}
\,\gamma^{\kern 0.5pt aq}_b\,\left(\!
\nabla_{\!q}\psi^{b1111\beta}[i]\,-
\vphantom{\sum^3_{\theta=1}}\right.\\
\vspace{2ex}
&\left.-\,\frac{2\,i\,e}{3\,\hbar\,c}\,A_q\,\psi^{b1111\beta}[i]
-\frac{i\,e\,g_3}{\hbar\,c}\,\sum^3_{\theta=1}
\uuuCA^\alpha_{q\kern 0.4pt\theta}
\,\psi^{b1111\theta}[i]\right)\,dV.
\endaligned\quad
\mytag{7.14}
$$ 
Looking at \mythetag{7.14}, we see that all upper level quarks, 
\pagebreak i\.\,e\. an up-quark, a charm-quark, and a top-quark 
(see table \mythetag{1.2}) have the same positive electric charge
$$
\hskip -2em
Q=+\frac{2}{3}\,e.
\mytag{7.15}
$$
Similarly, looking at \mythetag{7.13}, we find that all lower level 
quarks, i\.\,e\. a down-quark, a strange-quark, and a bottom-quark 
have the same negative electric charge
$$
\hskip -2em
Q=-\frac{1}{3}\,e.
\mytag{7.16}
$$
Compare \mythetag{7.15} and \mythetag{7.16} with \mythetag{6.15} 
in the case of charged leptons $e$, $\mu$, and $\tau$.\par
     The term $\Cal L_{83}$ in \mythetag{7.12} is a purely potential
term. It describes the interaction of lower level quarks with $Z$-bosons.
From \mythetag{7.10} and \mythetag{7.11} we derive
$$
\gathered
\Cal L_{83}=-\frac{e}{c}\int\!\!\sum_{i=d,s,b}\,
\sum^4_{a=1}\sum^4_{\bar a=1}\sum^4_{b=1}\sum^4_{c=1}\sum^3_{q=0}
\sum^3_{\beta=1}\sum^3_{\bar\beta=1}
\uD^{11}\,\uD^{11}\,\uuuD_{\beta\bar\beta}\,\times\\
\vspace{2ex}
\times\,\overline{\psi^{\raise 0.6pt \hbox{$\ssize\bar a
\bar\beta$}}_{11}[i]}
\,D_{a\bar a}\,\gamma^{\kern 0.5pt aq}_c\,Z_q\,
\frac{6\,(g_1)^2\,\circH^c_b-(3\,(g_1)^2+(g_2)^2)\,\bulletH^c_b}
{\sqrt{(g_2)^2+(3\,g_1)^2}}\ \psi^{\,b\beta}_{11}[i]\,dV.
\endgathered\quad
\mytag{7.17}
$$
The term $\Cal L_{84}$ in \mythetag{7.12} is also a purely potential
term. It describes the interaction of lower upper quarks with $Z$-bosons.
From \mythetag{7.9} and \mythetag{7.11} we derive
$$
\gathered
\Cal L_{84}=-\frac{e}{c}\int\!\!\sum_{i=d,s,b}\,
\sum^4_{a=1}\sum^4_{\bar a=1}\sum^4_{b=1}\sum^4_{c=1}\sum^3_{q=0}
\sum^3_{\beta=1}\sum^3_{\bar\beta=1}\uD^{11}\,\uD^{11}\,\times\\
\vspace{2ex}
\times\,\uD^{11}\,\uD^{11}\,\uuuD_{\beta\bar\beta}
\ \overline{\psi^{\raise 0.6pt \hbox{$\ssize\bar a1111
\bar\beta$}}[i]}\ D_{a\bar a}\,\gamma^{\kern 0.5pt aq}_c\,Z_q
\,\times\\
\vspace{2ex}
\times\,\frac{-12\,(g_1)^2\,\circH^c_b-(3\,(g_1)^2-(g_2)^2)
\,\bulletH^c_b}{\sqrt{(g_2)^2+(3\,g_1)^2}}
\ \psi^{\,b1111\beta}[i]\,dV.
\endgathered\quad
\mytag{7.18}
$$
The term $\Cal L_{85}$ in \mythetag{7.12} is a purely potential term
describing the interaction of upper and lower level quarks with 
$W$-bosons. From \mythetag{7.11} we derive
$$
\gathered
\Cal L_{85}=\frac{e}{c}\int\sum^3_{i=1}\,
\sum^4_{a=1}\sum^4_{\bar a=1}\sum^4_{b=1}\sum^4_{c=1}\sum^3_{q=0}
\sum^3_{\beta=1}\sum^3_{\bar\beta=1}
\uD^{11}\,\uD^{11}\,\uuuD_{\beta\bar\beta}\,\times\\
\vspace{2ex}
\times\,\overline{\psi^{\raise 0.6pt \hbox{$\ssize
\kern 0.5pt\bar a\bar\beta$}}_{11}[i]}
\,D_{a\bar a}\,\gamma^{\kern 0.5pt aq}_c\,g_2\,
W^{\sssize -}_{q111111}\,\bulletH^c_b
\ \psi^{\kern 0.5pt b1111\beta}[i]\,dV\,+\\
\vspace{2ex}
+\frac{e}{c}\int\sum^3_{i=1}\,\sum^4_{a=1}
\sum^4_{\bar a=1}\sum^4_{b=1}\sum^4_{c=1}\sum^3_{q=0}
\sum^3_{\beta=1}\sum^3_{\bar\beta=1}\uD^{11}\,\uD^{11}
\,\uD^{11}\,\times\\
\vspace{2ex}
\times\,\uD^{11}\,\uuuD_{\beta\bar\beta}
\ \overline{\psi^{\raise 0.6pt \hbox{$\ssize\kern 0.5pt
\bar a1111\bar\beta$}}[i]}\,D_{a\bar a}
\,\gamma^{\kern 0.5pt aq}_c\,g_2\,W^{{\sssize +}111111}_q
\,\bulletH^c_b\ \psi^{\kern 0.5pt b\beta}_{11}[i]\,dV.
\endgathered\quad
\mytag{7.19}
$$
The quark terms \mythetag{7.17}, \mythetag{7.18}, and \mythetag{7.19}
in the action integral are analogous to the terms \mythetag{6.17},
\mythetag{6.18}, and \mythetag{6.19} in the case of leptons. Note
that the terms responsible for interaction of quarks and gluons are
comprised within the kinetic terms \mythetag{7.13} and \mythetag{7.14}.
This means that the color symmetry $\MatGrSU(3)$ is not broken.\par
    Now let's proceed with the integrals \mythetag{7.2}. They are 
responsible for the masses of quarks. Substituting \mythetag{7.3}
into \mythetag{7.2}, we derive
$$
\hskip -2em
\Cal L_9=\Cal L_{91}+\Cal L_{92}.
\mytag{7.20}
$$
The first term $\Cal L_{91}$ in the sum \mythetag{7.20} is written
as follows:
$$
\gathered
\Cal L_{91}=-\sum^3_{i=1}\sum^3_{j=1}\frac{v}{\sqrt{2}}\int
\sum^4_{a=1}\sum^4_{\bar a=1}\sum^4_{b=1}\sum^3_{\beta=1}
\sum^3_{\bar\beta=1}
\uD_{11}\,\uD_{11}\,\uD_{11}\,\uD_{11}\,\times\\
\vspace{2ex}
\times\,\,\uuuD_{\beta\bar\beta}\,D_{a\bar a}
\left(h_1[ij]\,\bulletH^a_b+\overline{h_1[j\kern 0.5pt i]
\vphantom{A^A}}\,\circH^a_b\right)\overline{\psi^{\raise 0.6pt
\hbox{$\ssize\bar a1111\bar\beta$}}[i]}
\ \psi^{\kern 0.5pt b\kern 0.2pt 1111\beta}[j]\,dV\,-\\
\vspace{2ex}
-\sum^3_{i=1}\sum^3_{j=1}\frac{v}{\sqrt{2}}\int\sum^4_{a=1}
\sum^4_{\bar a=1}\sum^4_{b=1}\sum^3_{\beta=1}\sum^3_{\bar\beta=1}
\uD^{11}\,\uD^{11}\,\times\\
\vspace{2ex}
\times\,\,\uuuD_{\beta\bar\beta}\,D_{a\bar a}
\left(h_2[ij]\,\bulletH^a_b+\overline{h_2[j\kern 0.5pt i]
\vphantom{A^A}}\,\circH^a_b\right)
\overline{\psi^{\raise 0.6pt \hbox{$\ssize\kern 0.5pt\bar a
\bar\beta$}}_{11}[i]}
\ \psi^{\kern 0.5pt b\kern 0.2pt\beta}_{11}[j]\,dV.
\endgathered\quad
\mytag{7.21}
$$
The second term $\Cal L_{92}$ in \mythetag{7.20} is responcible for
the interaction of quarks with the real scalar Higgs field $\chi$. 
It is very similar to \mythetag{7.21}:
$$
\gathered
\Cal L_{92}=-\sum^3_{i=1}\sum^3_{j=1}\frac{1}{\sqrt{2}}\int
\sum^4_{a=1}\sum^4_{\bar a=1}\sum^4_{b=1}\sum^3_{\beta=1}
\sum^3_{\bar\beta=1}
\uD_{11}\,\uD_{11}\,\uD_{11}\,\uD_{11}\,\times\\
\vspace{2ex}
\times\,\,\uuuD_{\beta\bar\beta}\,D_{a\bar a}\,\chi
\left(h_1[ij]\,\bulletH^a_b+\overline{h_1[j\kern 0.5pt i]
\vphantom{A^A}}\,\circH^a_b\right)\overline{\psi^{\raise 0.6pt
\hbox{$\ssize\bar a1111\bar\beta$}}[i]}
\ \psi^{\kern 0.5pt b\kern 0.2pt 1111\beta}[j]\,dV\,-\\
\vspace{2ex}
-\sum^3_{i=1}\sum^3_{j=1}\frac{1}{\sqrt{2}}\int\sum^4_{a=1}
\sum^4_{\bar a=1}\sum^4_{b=1}\sum^3_{\beta=1}\sum^3_{\bar\beta=1}
\uD^{11}\,\uD^{11}\,\times\\
\vspace{2ex}
\times\,\,\uuuD_{\beta\bar\beta}\,D_{a\bar a}\,\chi
\left(h_2[ij]\,\bulletH^a_b+\overline{h_2[j\kern 0.5pt i]
\vphantom{A^A}}\,\circH^a_b\right)
\overline{\psi^{\raise 0.6pt \hbox{$\ssize\kern 0.5pt\bar a
\bar\beta$}}_{11}[i]}
\ \psi^{\kern 0.5pt b\kern 0.2pt\beta}_{11}[j]\,dV.
\endgathered\quad
\mytag{7.22}
$$
The term \mythetag{7.21} is a mass term for quarks. However, in
general case, using it, one cannot prescribe masses to individual
quarks. Let's consider a special case, where the diagonal elements
of the coupling constants matrices are real constants:
$$
\xalignat 2
&\hskip -2em
h_1[i\kern 0.5pt i]=\overline{h_1[i\kern 0.5pt i]\vphantom{A^A}},
&&h_2[i\kern 0.5pt i]=\overline{h_1[i\kern 0.5pt i]\vphantom{A^A}}.
\quad
\mytag{7.23}
\endxalignat
$$
In this special case, where the equalities \mythetag{7.23} are fulfilled,
we can write
$$
\hskip -2em
\Cal L_{91}=\Cal L_{91}[diag]+\Cal L_{91}[not\,diag].
\mytag{7.24}
$$
The diagonal term $\Cal L_{91}[diag]$ in the expansion \mythetag{7.24}
is written as follows:
$$
\gathered
\Cal L_{91}[diag]=-\sum^3_{i=1}\frac{h_1[i\kern 0.5pt
i]\,v}{\sqrt{2}}\int\sum^4_{a=1}\sum^4_{\bar a=1}\sum^3_{\beta=1}
\sum^3_{\bar\beta=1}\uD_{11}\,\uD_{11}\,\uD_{11}\,\uD_{11}\,\times\\
\vspace{2ex}
\times\,\uuuD_{\beta\bar\beta}\,D_{a\bar a}
\,\overline{\psi^{\raise 0.6pt \hbox{$\ssize\bar a1111\bar\beta$}}[i]}
\ \psi^{\kern 0.5pt a\kern 0.2pt 1111\beta}[i]\,dV
-\sum^3_{i=1}\frac{h_1[i\kern 0.5pt
i]\,v}{\sqrt{2}}\,\times\\
\vspace{2ex}
\times\int\sum^4_{a=1}\sum^4_{\bar a=1}
\sum^3_{\beta=1}\sum^3_{\bar\beta=1}
\uD^{11}\,\uD^{11}\,\uuuD_{\beta\bar\beta}\,D_{a\bar a}
\,\overline{\psi^{\raise 0.6pt \hbox{$\ssize\kern 0.5pt\bar a
\bar\beta$}}_{11}[i]}
\ \psi^{\kern 0.5pt b\kern 0.2pt\beta}_{11}[i]\,dV.
\endgathered\qquad
\mytag{7.25}
$$
Looking at \mythetag{7.25}, one can prescribe the following masses
to individual quarks:
$$
\xalignat 3
&\hskip -2em
m_u=\frac{h_1[11]\,v}{\sqrt{2}\,c}\,,
&&m_c=\frac{h_1[22]\,v}{\sqrt{2}\,c}\,,
&&m_t=\frac{h_1[33]\,v}{\sqrt{2}\,c}\,,
\quad\\
\vspace{-1ex}
&&&\mytag{7.26}\\
\vspace{-1ex}
&\hskip -2em
m_d=\frac{h_2[11]\,v}{\sqrt{2}\,c}\,,
&&m_s=\frac{h_2[22]\,v}{\sqrt{2}\,c}\,,
&&m_b=\frac{h_2[33]\,v}{\sqrt{2}\,c}\,.
\quad
\endxalignat
$$
The formulas \mythetag{7.26} are similar to the formulas \mythetag{6.23}
for charged leptons. Note that we can impose a more restrictive condition
for the coupling constants than that of \mythetag{7.23}, e\.\,g\. we can
require them to form Hermitian matrices:
$$
\xalignat 2
&\hskip -2em
h_1[ij]=\overline{h_1[j\kern 0.5pt i]\vphantom{A^A}},
&&h_2[ij]=\overline{h_1[j\kern 0.5pt i]\vphantom{A^A}}.
\quad
\mytag{7.27}
\endxalignat
$$
In this case the action integrals \mythetag{7.21} and
\mythetag{7.22} simplify to
$$
\gather
\hskip -2em
\Cal L_{91}=-\sum^3_{i=1}\sum^3_{j=1}\frac{h_1[ij]\,v}{\sqrt{2}}\int
\sum^4_{a=1}\sum^4_{\bar a=1}\sum^4_{b=1}\sum^3_{\beta=1}
\sum^3_{\bar\beta=1}
\uD_{11}\,\uD_{11}\,\uD_{11}\,\uD_{11}\,\times\quad\\
\vspace{1ex}
\hskip -2em
\times\,\,\uuuD_{\beta\bar\beta}\,D_{a\bar a}
\ \overline{\psi^{\raise 0.6pt
\hbox{$\ssize\bar a1111\bar\beta$}}[i]}
\ \psi^{\kern 0.5pt b\kern 0.2pt 1111\beta}[j]\,dV
-\sum^3_{i=1}\sum^3_{j=1}\frac{h_1[ij]\,v}{\sqrt{2}}\,\times\quad
\mytag{7.28}\\
\vspace{1ex}
\hskip -2em
\times\int\sum^4_{a=1}
\sum^4_{\bar a=1}\sum^4_{b=1}\sum^3_{\beta=1}\sum^3_{\bar\beta=1}
\uD^{11}\,\uD^{11}\,\uuuD_{\beta\bar\beta}\,D_{a\bar a}
\ \overline{\psi^{\raise 0.6pt \hbox{$\ssize\kern 0.5pt\bar a
\bar\beta$}}_{11}[i]}
\ \psi^{\kern 0.5pt b\kern 0.2pt\beta}_{11}[j]\,dV,\quad\\
\vspace{2ex}
\hskip -2em
\Cal L_{92}=-\sum^3_{i=1}\sum^3_{j=1}\frac{h_1[ij]}{\sqrt{2}}\int
\sum^4_{a=1}\sum^4_{\bar a=1}\sum^4_{b=1}\sum^3_{\beta=1}
\sum^3_{\bar\beta=1}
\uD_{11}\,\uD_{11}\,\uD_{11}\,\uD_{11}\,\times\kern -20pt\quad\\
\vspace{1ex}
\hskip -2em
\times\,\,\uuuD_{\beta\bar\beta}\,D_{a\bar a}
\,\chi\ \overline{\psi^{\raise 0.6pt
\hbox{$\ssize\bar a1111\bar\beta$}}[i]}
\ \psi^{\kern 0.5pt b\kern 0.2pt 1111\beta}[j]\,dV
-\sum^3_{i=1}\sum^3_{j=1}\frac{h_1[ij]}{\sqrt{2}}\,\times\quad
\mytag{7.29}\\
\vspace{1ex}
\hskip -2em
\times\int\sum^4_{a=1}
\sum^4_{\bar a=1}\sum^4_{b=1}\sum^3_{\beta=1}\sum^3_{\bar\beta=1}
\uD^{11}\,\uD^{11}\,\uuuD_{\beta\bar\beta}\,D_{a\bar a}\,\chi
\ \overline{\psi^{\raise 0.6pt \hbox{$\ssize\kern 0.5pt\bar a
\bar\beta$}}_{11}[i]}
\ \psi^{\kern 0.5pt b\kern 0.2pt\beta}_{11}[j]\,dV.\quad
\endgather
$$
Since \mythetag{7.27} implies \mythetag{7.23}, the formulas 
\mythetag{7.26} for quark masses are valid in this case too.
Unlike \mythetag{7.21} and \mythetag{7.22}, the action integrals
\mythetag{7.28} and \mythetag{7.29} preserve the chiral-to-antichiral
symmetry. However, this makes no difference for the Standard Model
in whole since there are many other terms in the total action 
integral that break this symmetry.
\head
8. Conclusion.
\endhead
     The main purpose of this paper is to explain the Standard
Model of elementary particles in a little bit non-standard way
different from that traditionally used in physical literature. 
In two previous papers \mycite{2} and \mycite{6} three special
complex vector bundles over the space-time manifold $M$ were 
introduced and studied. These bundles provide a geometric 
background for describing the Standard Model in the case of a
non-flat space-time manifold $M$, i\.\,e\. in the presence of
a gravitation field. The actual description of the Standard Model
in terms of these three bundles is given in the present paper. 


\Refs
\ref\myrefno{1}\by Sharipov~R.~A.\paper A note on Dirac spinors 
in a non-flat space-time of general relativity\publ e-print 
\myhref{http://arXiv.org/abs/math/0601262/}{math.DG/0601262} 
in Electronic Archive \myEarXivlink
\endref
\ref\myrefno{2}\by Sharipov~R.~A.\paper The electro-weak and color 
bundles for the Standard Model in a gravitation field\publ e-print 
\myhref{http://arXiv.org/abs/math/0603611/}{math.DG/0603611} 
in Electronic Archive \myEarXivlink
\endref
\ref\myrefno{3}\by Rubakov~V.~A.\book Classical gauge fields
\publ Editorial URSS\publaddr Moscow\yr 1999
\endref
\ref\myrefno{4}\by Cianfrani~F., Montani~G.\paper Geometrization 
of the electro-weak model bosonic component\publ\linebreak e-print
\myhref{http://arXiv.org/abs/gr-qc/0601052/}{gr-qc/0601052} 
in Electronic Archive \myEarXivlink
\endref
\ref\myrefno{5}\by Sharipov~R.~A.\paper On the Dirac equation in 
a gravitation field and the secondary quantization\publ e-print 
\myhref{http://uk.arXiv.org/abs/math/0603367/}{math.DG/0603367} 
in Electronic Archive \myEarXivlink
\endref
\ref\myrefno{6}\by Sharipov~R.~A.\paper A note on connections 
of the Standard Model in a gravitation field\publ e-print 
\myhref{http://uk.arXiv.org/abs/math/0604145/}{math.DG/0604145} 
in Electronic Archive \myEarXivlink
\endref
\ref\myrefno{7}\by Sharipov~R.~A.\paper A note on metric connections 
for chiral and Dirac spinors\publ e-print 
\myhref{http://arXiv.org/abs/math/0602359/}{math.DG}
\myhref{http://arXiv.org/abs/math/0602359/}{/0602359}
in Electronic Archive \myEarXivlink
\endref
\ref\myrefno{8}\by Sharipov~R.~A.\book Classical electrodynamics and
theory of relativity\publ Bashkir State University\publaddr Ufa\yr 1997
\moreref see also
\myhref{http://arXiv.org/abs/physics/0311011}{physics/0311011}
in Electronic Archive \myEarXivlink\ and 
\myhref{http://www.geocities.com/r-sharipov/r4-b5.htm}
{r-sharipov/r4-} \myhref{http://www.geocities.com/r-sharipov/r4-b5.htm}
{b5.htm} in \myGeoCities
\endref
\endRefs
\enddocument
\end